\documentclass[reqno,11pt]{amsart}
\usepackage{amsmath,amsthm,amscd,amsfonts,amssymb,hyperref,floatrow}
\usepackage{graphicx, color,dsfont,xcolor}
\usepackage{mathtools}

\numberwithin{equation}{section}

\newtheorem{hypothesis}{Hypothesis}

\DeclareMathSymbol{\leqslant}{\mathalpha}{AMSa}{"36} % nicer `smaller or equal'
\DeclareMathSymbol{\geqslant}{\mathalpha}{AMSa}{"3E} % nicer `larger or equal'
\DeclareMathSymbol{\eset}{\mathalpha}{AMSb}{"3F}     % nicer `emptyset'
\renewcommand{\leq}{\;\leqslant\;}                   % redef. of < or =
\def \C{ \mathbb  C }
\newcommand{\R}{\mathbb{R}}
\newcommand{\Z}{\mathbb{Z}}
\newcommand{\N}{\mathbb{N}}

\newcommand{\ZZ}{\mathbb{Z}}

\def \E{ \mathbb E  }

\begin{document}
\title{The eigenvectors of Gaussian matrices with an external source}

\author{Romain Allez \and Jo\"el Bun \and Jean-Philippe Bouchaud}
\address{Weierstrass Institute, Mohrenstr. 39, 10117 Berlin, Germany.}
\address{Capital Fund Management, 23-25, rue de l'Universit\'e, 75007 Paris}
\address{CNRS, LPTMS, Batiment 100, Universit\'e d'Orsay, 91405 Orsay Cedex, France}
\address{Leonard de Vinci P{\^ o}le Universitaire, Finance Lab, 92916 Paris La D{\'e}fense, France}
\email{romain.allez@gmail.com} 
\email{joel.bun@gmail.com} 
\email{jean-philippe.bouchaud@cfm.fr}
\date{\today}
\maketitle

\begin{abstract}
We consider a diffusive matrix process $(X_t)_{t\ge 0}$  
defined as $X_t:=A+H_t$ where $A$ is a given deterministic Hermitian 
matrix and $(H_t)_{t\ge 0}$ is a Hermitian Brownian motion. 
The matrix $A$ is the ``external source'' that one would like to estimate
from the noisy observation $X_t$ at some time $t>0$. 
We investigate the relationship between the non-perturbed eigenvectors of the matrix $A$
and the perturbed eigenstates at some time $t$ for the three 
relevant scaling relations between the time $t$ and the dimension $N$ of the matrix $X_t$. 
We determine the asymptotic (mean-squared) projections of any given 
non-perturbed eigenvector $|\psi_j^0\rangle$, 
associated to an eigenvalue $a_j$ of $A$ which may lie inside the bulk of the spectrum or be isolated (spike)
from the other eigenvalues, on the orthonormal basis of the perturbed eigenvectors $|\psi_i^t\rangle,i\neq j$. 
We derive a Burgers type evolution equation for the local resolvent $(z-X_t)_{ii}^{-1}$, describing 
the evolution of the local density of a given initial state $|\psi_j ^0\rangle$. We are able to solve this equation explicitly 
in the large $N$ limit, for any initial matrix $A$.
In the case of one isolated eigenvector $|\psi_j^0\rangle$, 
we prove a central limit Theorem for the overlap $\langle \psi_j^0|\psi_j^t\rangle$. 
When properly centered and rescaled by a factor $\sqrt{N}$, this overlap converges in law towards a centered 
Gaussian distribution with an explicit variance depending on $t$.  
Our method is based on analyzing the eigenvector flow under the Dyson Brownian motion. 
 \end{abstract}

\vspace{1cm}

\section{Introduction}
We consider a diffusion process $(X_t)_{t\ge 0}$ in the space of $N\times N$ symmetric or Hermitian matrices
starting from a given deterministic matrix $A$ and evolving with time according to a Hermitian Brownian motion. 
The matrix $X_t$ at time $t$ is  
\begin{align}\label{def-X}
X_t:= A + H_t
\end{align}
where $(H_t)_{t\ge 0}$ is a Hermitian Brownian motion, i.e. a diffusive matrix process such that  $H_0=0$
and whose entries $\{H_{t}(ij), i \le j\}$ are given by
\begin{align}\label{hermitian-bm}
H_t(ij):= \frac{1}{\sqrt{N}} B_t(ij) \quad \mbox{  if  } \quad  i \neq j,  \quad H_t(ii):= \frac{\sqrt{2}}{\sqrt{N}} B_t(ii)
\end{align}
where the $B_t(ij), i \le j$ are independent and identically distributed real or complex (real if $i=j$) Brownian motions.  

The matrix $A$ is the external source of the title and should be seen as a signal that one would like to estimate 
from the observation of the noisy matrix $X_t$.   

In this paper, we are interested in the effect of the addition of the noisy perturbation matrix $H_t$ in the limit of large dimension $N\to +\infty$. 
More precisely, we investigate the relationship between the eigenvectors of the perturbed matrix $X_t$ with those of the 
initial matrix $A$ for some given $t>0$, possibly scaling with the dimension $N$ of the matrices. 
We will see that there are in fact several regimes to study depending on the scaling relation  
between $t$ and $N$ and on the positions of the associated eigenvalues which may be 
isolated or inside the continuous part (bulk) of the spectrum.   

The evolution as $t$ grows of the eigenvalues $\lambda_1(t) \ge \lambda_2(t) \ge \dots \ge \lambda_N(t)$ 
of the symmetric  matrix $X_t$  has been investigated in tremendous details
in random matrix theory (see \cite[section 4.3]{agz} for a review). 
It was first shown by Dyson \cite{dyson} in 1962  that the eigenvalues of the matrix $X_t$ 
evolve according to the Dyson Brownian motion which describes the evolution of  
 $N$ positively charged particles (Coulomb gas) subject to electrostatic repulsion and 
 to independent thermal noises. 
 The dynamics of the Dyson Brownian motion were studied in many details for different purposes. 
 The most striking applications of the Dyson Brownian motion 
 are perhaps the proofs of the universality conjectures for Wigner matrices (see e.g. \cite{universality,bourgade-univ} 
 and references therein). The Dyson Brownian motion was also used in theoretical physics as a model 
 to study disordered metals and chaotic billiards \cite{beenakker} (see also \cite{forrester}). In this context, the authors 
 compute the correlations between the positions of the eigenvalues in the bulk 
 at a given time $s$ with those at a later time $t>s$. 
The asymptotic correlation functions are described in terms of the extended Hermite kernel. 
The correlations between the positions of the eigenvalues near the edge of the spectrum at different times 
were later computed in \cite{macedo} in terms of the extended Airy kernel.  

The study of the associated eigenvectors denoted respectively by $|\psi_1^t\rangle, |\psi_2^t\rangle,\dots, |\psi_N^t\rangle$
is comparatively much poorer. A few authors were interested in some aspects 
of eigenvector fluctuations (see e.g. in \cite{donati,florent,bourgade-yau,tao} on the statistics of Haar matrices, 
\cite{sandrine,paul,benaych} 
for eigenvectors of covariance matrices and \cite{vectors} for applications in finance)
but yet very little is known about the cross correlation of the eigenvectors at different times $s$ and $t>s$. 
It is a natural question to extend the results known for the eigenvalues \cite{beenakker,macedo} by investigating 
the relation between the eigenvectors of the matrix $X_s$ with those 
at a later time $t>s$ (with possibly $s=0$). 
This question was initiated in \cite{wilkinson1} and recently reconsidered in \cite{vectors2} where we investigate the projections of a given eigenvector 
$|\psi_i ^0\rangle$ at time $0$ on the orthonormal basis of the perturbed eigenvectors at time $t$. 
Specifically, we consider the case where the associated eigenvalue $\lambda_i(0)$ lies in the continuous part 
of the spectrum and use Stieltjes transform methods to compute the asymptotic (mean squared) 
projections of this vector 
on the orthonormal  basis at time $s=0$. In this paper, we use a powerful method 
based on analyzing the eigenvector flow under the Dyson Brownian motion. 
This method was used previously in \cite{bourgade-univ}. We obtain an autonomous 
equation \eqref{eq-overlaps} 
satisfied by the projections of a  given non-perturbed eigenvector $|\psi_j^0\rangle$ 
on the basis of the perturbed eigenstates $|\psi_i^t\rangle$ at time $t$. 
This approach permits us to analyze the relation between this non-perturbed deterministic and fixed vector $|\psi_j^0\rangle$ 
(to be estimated in applications)
with the perturbed (random) states, in a general setting
for the locations of the associated non-perturbed eigenvalue $\lambda_j (0)$ (inside or outside the bulk)
and for the different relevant scaling regimes between the time $t$ and the dimension $N$.

In order to explain more precisely the general extent of our results, let us make an assumption on the spectrum of the 
initial matrix $A$. 

\begin{figure}\center
\includegraphics[scale=0.8]{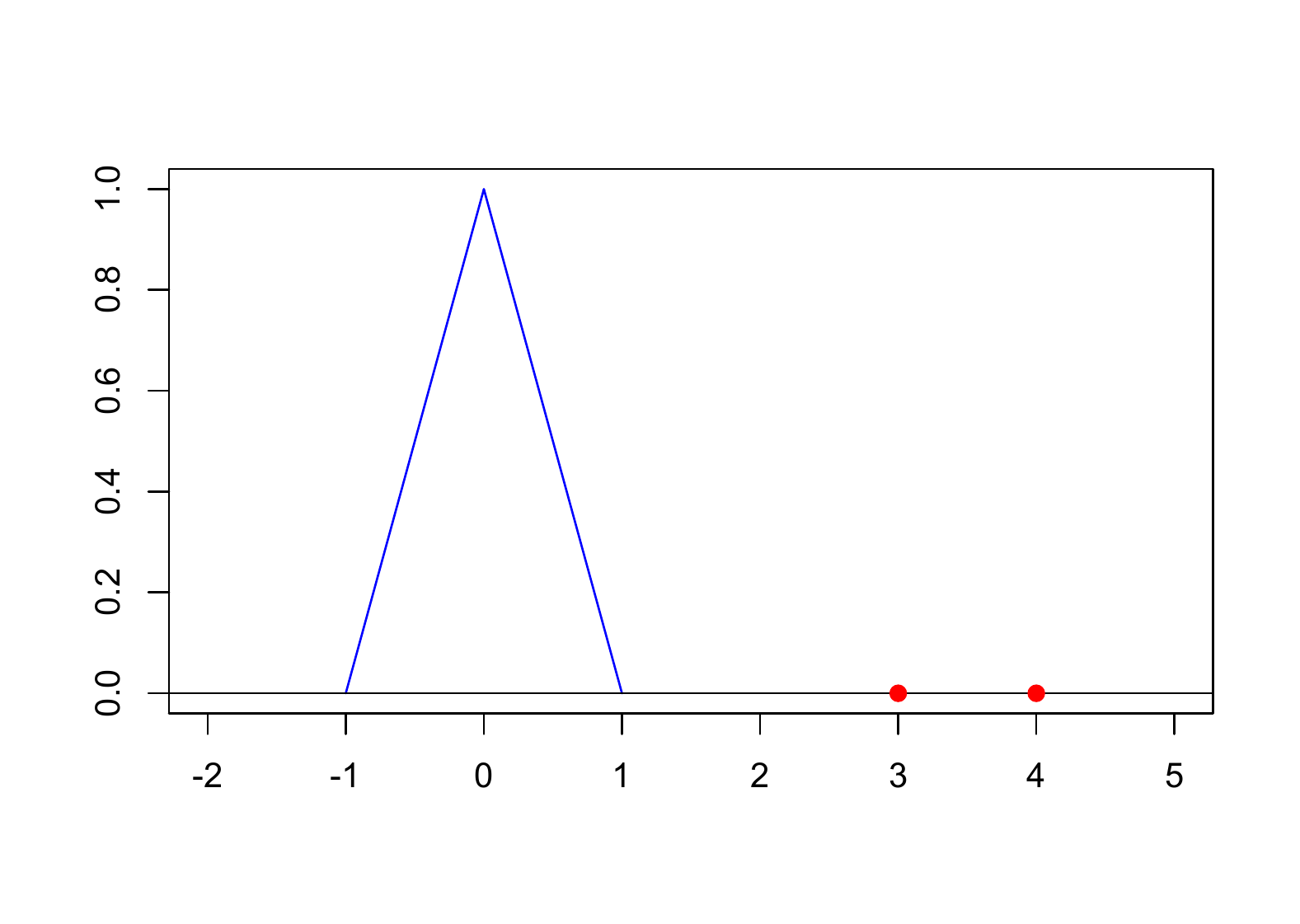}
\caption{Spectrum of a matrix $A$ satisfying our hypothesis with a continuous triangular density and  
$k=2$ spikes $a_1=4, a_2=3$. }\label{fig-triangle}
\end{figure}

\begin{hypothesis}\label{hypothesis}
We consider a sequence of symmetric matrices $A:=(A_N)_{N\in \N}$ such that $A_N$ has size $N\times N$. 
We denote by $a_1\ge a_2\ge \dots \ge a_N$ the eigenvalues of the matrix $A_N$ 
\footnote{To simplify notations, we do not use an additional superscript $N$.}
and by  $|\phi_1\rangle, |\phi_2 \rangle, \dots,|\phi_N\rangle$
the respective associated eigenvectors. 
We will work under the assumption that the spectrum of $A$ 
can be decomposed in the large $N$ limit into a continuous and a discrete part. 
More precisely, we will suppose that 
\begin{itemize}
\item The empirical eigenvalue density of the matrix $A$ converges weakly 
to some limiting compactly supported density $\rho_A(\lambda)$ when $N\to \infty$ i.e.
\begin{align}\label{conv-rho-A}
\frac{1}{N} \sum_{i=1}^N \delta_{a_i}(d\lambda) \to \rho_A(\lambda) \, d\lambda\,. 
\end{align}
The compactness of the support is not necessary for the results established below but 
we need this assumption in order to have some empty space left for additional isolated eigenvalues.  
We will also denote by $a(x)$ the smooth function such that for any $x\in [0,1]$,
\begin{align}\label{definition-a}
x = \int_{a(x)}^{+\infty} \rho_A(\lambda)\, d\lambda\,.  
\end{align}
If $x\in [0,1]$, $a(x)$ is the $x$-quantile of the probability distribution $\rho_A$. 
To avoid technicalities, we will even work under a stronger assumption than the convergence \eqref{conv-rho-A} 
related to the local repartition 
of the eigenvalues of $A$. We will suppose that the eigenvalues of $A$ are allocated smoothly according 
to the quantile of the limiting smooth \footnote{The probability density $\rho_A$ is supposed to be smooth 
at least on the interior of the compact support.} density $\rho_A$
i.e. such that for any $k=1,\dots,N,$ 
\begin{align*}
a_k = a(\frac{k}{N}).
\end{align*}
\item In addition to the compactly supported probability density $\rho_A$, we suppose that 
there is a fixed finite number $\ell$ (independent of $N$) 
of {\it isolated} eigenvalues (with multiplicity one). For simplicity, we will 
suppose without loss of generality that those ``spikes'' are the largest eigenvalues $a_1 > a_2 > \dots >a_\ell$ 
(see Fig. \ref{fig-triangle} for an illustration of a matrix $A$ satisfying hypothesis \ref{hypothesis}).
\end{itemize}  
\end{hypothesis}

Let us now dress a brief outline of the paper. In section \ref{eigenvalues-traj}, 
we review some known facts on the evolution of the eigenvalues of $X_t$ when the dimension $N$ is finite and then in the large $N$ limit. We explain in subsection \ref{subsection-eq-overlaps} the derivation of the autonomous equation 
\eqref{eq-overlaps} which describes the evolution of
 the projections of a  given non-perturbed eigenvector $|\psi_j^0\rangle$ 
on the basis of the perturbed eigenstates $|\psi_i^t\rangle$ as a function of time $t$.
This equation \eqref{eq-overlaps}  is  a particular case of the eigenvector moment flow discovered in \cite{bourgade-yau}  
and will be our main tool to analyze the eigenvector dynamics in the next sections. 
In section \ref{bulk}, we investigate the overlaps of a given vector $|\psi_j^0\rangle$ with the perturbed eigenstates $|\psi_i^t\rangle$ 
at time $t$
when the associated eigenvalue $\lambda_j(0)=a_j$ lies in the bulk of the spectrum of the matrix $A$. 
We distinguish three different regimes for the values of $t$ as a function of $N$. The first regime is perturbative and corresponds 
to microscopic values of $t\ll 1/N$. The second mesoscopic regime, corresponding to values of $t$ proportional to $1/N$, 
establishes a smooth crossover between the first and the third regimes. 
A transition occurs in this mesoscopic regime: if $t\propto\tau/N$, the non-perturbed 
eigenvector $|\psi_j^0\rangle$ associated to $\lambda_j(0)=a_j$ is localized in the basis of the perturbed states 
$|\psi_i^t\rangle$ for finite values of $\tau$ while it becomes delocalized when $\tau \to +\infty$. 
The third regime corresponds to macroscopic values of $t$ which do not depend on $N$: $t$ is fixed and $N\to \infty$. 
To compute the limiting overlaps in this regime, we derive a new Burgers 
type evolution equation for the local resolvent $(z-X_t)_{ii}^{-1}$ of the matrix $X_t$. 
This equation characterizes the evolution of the local density of any given initial state $|\psi_j^0\rangle$
in the basis of the perturbed eigenvectors $|\psi_i^t\rangle, i=1,\dots,N$.
We are able to solve this local Burgers equation explicitly in the scaling limit, providing a nice explicit solution \eqref{sol-non-perturb}
to the overlap equation \eqref{overlap-eq-non-perturb}. We also check the perfect matching between the two Formulas
\eqref{cauchy} and \eqref{sol-non-perturb} which respectively describe the limiting overlaps 
in the mesoscopic and macroscopic regimes, at the frontier between those two regimes.
We provide many other details on the overlap statistics when $N\to \infty$ 
in the second and third regimes in subsections \ref{crossover} and \ref{third-regime}.
We consider the case where $A$ satisfies the general  hypothesis \ref{hypothesis} and also some special cases,
which are analyzed explicitly.  
We also give a few interesting open questions for future research. In section \ref{isolated}, 
we fully analyze the overlap $\langle \psi_1^0 | \psi_1^t\rangle$ 
between the non-perturbed eigenvector $|\psi_1^0\rangle$ associated to the largest spike 
of the matrix $A$ and the perturbed eigenvector $|\psi_1^t\rangle$ associated to the largest eigenvalue $\lambda_1(t)$ 
of the matrix $X_t$ \footnote{Our study is also valid for  the other eigenvectors $|\psi_2\rangle, |\psi_3\rangle$ 
associated to the next spikes.}. More precisely, we rigorously prove that 
the overlap $\langle \psi_1^0 |\psi_1^t\rangle$ properly centered and rescaled by a factor $\sqrt{N}$ 
converges in law when $N\to \infty$
towards a Gaussian distribution with an explicit variance depending on $t$ and on the eigenvalues trajectories. 
We are also able to characterize the variance of the 
transverse components of the non-perturbed 
eigenvector $|\psi_1^0\rangle$ carried by the perturbed eigenvectors $|\psi_i^t\rangle$ associated to the eigenvalues in the 
bulk of the spectrum of $X_t$. We finally consider the special case where $A$ has a small rank compared to $N$ (factor model)
in subsection \ref{factor-model-vectors}. 
In this case, we can compute the overlaps between the perturbed and non-perturbed states analytically. 
In the last section \ref{covariance}, 
we explain how our ideas could be used in the context of covariance matrices to compute the 
asymptotic overlaps between the sample and population eigenvectors when the dimension $N$ is very large.  

{\bf Acknowledgments}
We are grateful to Paul Bourgade, Antoine Dahlqvist, Laure Dumaz and Marc Potters for useful comments and discussions. 

RA received funding from the European Research Council under the European
Union's Seventh Framework Programme (FP7/2007-2013) / ERC grant agreement nr. 258237 and thanks the Statslab in DPMMS, Cambridge for its hospitality.

\section{Eigenvalues and eigenvectors trajectories}\label{eigenvalues-traj}

\subsection{Eigenvalues and eigenvectors diffusion processes}
It is well known  \cite{agz} that the eigenvalues $\lambda_1(t) \ge \lambda_2(t) \ge \dots \ge \lambda_N(t)$ 
of the matrix $X_t$ 
evolve according to the Dyson Brownian motion
\begin{align}\label{dyson-bm}
d\lambda_i(t) = \sqrt{\frac{2}{\beta N}} \,  dB_i(t) + \frac{1}{N} \sum_{j \neq i} \frac{dt}{\lambda_i(t) - \lambda_j(t)}\,, \quad
i=1,\dots, N\,,
\end{align}
where the $B_i$ are independent real Brownian motions, and satisfy the initial conditions 
\begin{align*}
\lambda_i(0)= a_i \,, \quad i=1,\dots, N\,. 
\end{align*} 
The eigenvalues of $X_t$ may be seen as 
positively charged particles in a one-dimensional Coulomb gas with electrostatic repulsion between them 
and subject to a thermal noise $dB_i(t)$.  

Conditionally on the eigenvalues paths, the trajectories of the associated eigenvectors
$ |\psi_1^t\rangle,$ $|\psi_2^t \rangle, \dots,$$|\psi_N^t\rangle$ can be realized continuously as a function of $t$. 
(We use Dirac's quantum mechanics bra-ket notation throughout this paper).
This eigenvector flow was first exhibited in \cite{bru} for Wishart processes.
Those continuous paths are determined using standard perturbation theory or stochastic analysis 
tools (see again \cite{agz} or \cite{alice}): in our case, we have, for all $i=1,\dots,N,$ 
\begin{align}\label{evolution-psi_i}
d|\psi_i^t\rangle& = - \frac{1}{2N}\sum_{j \neq i} \frac{dt}{(\lambda_i(t)-\lambda_j(t))^2} |\psi_i^t\rangle + 
\frac{1}{\sqrt{N}} \sum_{j \neq i } \frac{dw_{ij}(t)}{\lambda_i(t) - \lambda_j(t)} |\psi_j^t \rangle \,, \\
&\mbox{with } \quad  |\psi_i^0\rangle = |\phi_i\rangle \,,
\end{align} 
where the family of independent (up to symmetry) real Brownian motions $\{w_{ij}: i \neq j \}$ is independent 
of the eigenvalues trajectories (i.e. independent of the driving Brownian motions $B_i$ in \eqref{dyson-bm}).  
We can therefore freeze the eigenvalues trajectories and then, conditionally on this eigenvalues path, study the eigenvectors evolution. 
The eigenvector process can thus be regarded as a diffusion process in a random environment which depends on the realized trajectories of the eigenvalues. 
This is an important fact that will be used several times throughout this paper. Most of the results derived in this paper concern the large dimensional statistics
of the eigenvectors and hold almost surely with respect to the eigenvalues trajectories.  

The evolution equation \eqref{evolution-psi_i} for the $i$-th eigenvector contains two orthogonal terms. The first term, 
collinear to $|\psi_i^t\rangle$, pulls back $|\psi_i^t\rangle$ towards $0$ 
in such a way that the eigenvectors remain normalized $\langle \psi_i ^t|\psi_i^t\rangle=1$.  
The randomness comes in the second interaction and transverse term. 
We see that the $i$-th eigenvector $|\psi_i^t\rangle$ trades more information with the eigenvectors $|\psi_j^t\rangle, j\neq i$ 
that are associated to the closest neighboring eigenvalues $\lambda_j(t) \sim \lambda_i(t)$. 
If the neighboring eigenvalues $\lambda_j(t)$ are very close to $\lambda_i(t)$ 
(typically at a distance of order $1/N$ in the continuous part of the spectrum for large $N$), we shall see that 
this singular interaction leads to unstable (discontinuous) eigenstates trajectories 
with respect to time $t$, in the large $N$ limit (see below).

\subsection{Evolution of the mean squared overlaps at finite $N$}\label{subsection-eq-overlaps}
In order to quantify the relationship between the perturbed eigenstates at time $t$ 
and the eigenstates at the initial time, 
we consider the scalar products or {\it overlaps} $\langle \psi_i^t |\psi_j^0 \rangle$ for $i,j =1,\dots,N$.
Specifically, we investigate the mean square overlaps $\left [ \langle \psi_i^t |\psi_j^0 \rangle^2\right]$ where we use 
the notation $\left [ \dots \right] $ for the expectation over the Brownian motions $w_{ij}, i\neq j \in \{1,\dots,N\}$ 
which appear in the eigenvectors evolution equation \eqref{evolution-psi_i}. Recall that those 
Brownian motions are {\it independent} of the eigenvalues so that this conditioning does
not modify the law of the eigenvalue process. 
Note also that the variables $\left [ \langle \psi_i^t |\psi_j^0 \rangle^2\right], 1\le i,j\le N$ are still random, 
measurable with respect to the sigma field generated by the Brownian trajectories $\{(B_i(s)), 0\le s \le t,i=1,\cdots,N\}$.

A straightforward application of It\^o's formula permits us to find an evolution equation for the mean squared 
overlaps $\left[ \langle \psi_i ^t|\psi_j^0\rangle^2\right]$.
More precisely, if $j$ is fixed, one can show that the 
(mean squared) projections of the non-perturbed 
eigenstate $|\psi_j^0\rangle$  on the perturbed eigenvectors $|\psi_i^t\rangle, i =1,\dots,N$ satisfy an {\it autonomous} evolution equation. 
We will use the following short hand notation for the overlaps
%It is convenient to index the overlaps by the associated eigenvalues $\lambda_i(t)$ 
%and $a_j$ instead of the indices $i$ and $j$ by introducing the notation 
\begin{align}
u_{i|j}(t):= \left[ \langle \psi_i ^t|\psi_j^0\rangle^2\right] \,.
\end{align}   
 The overlap equation for the (mean squared) projections of the state $|\psi_j^0\rangle$ reads 
 \begin{align}\label{eq-overlaps}
\partial_t\, u_{i|j}(t) = \frac{1}{N} \sum_{k \neq i} \frac{u_{k|j}(t)- u_{i|j}(t)}{(\lambda_k(t)-\lambda_i(t))^2}\quad 
\mbox{ with } \quad u_{i|j}(0) = \delta_{ij}\,. 
\end{align} 
This evolution equation was discovered in 1995 by Wilkinson and Walker (see Eq. (4.7) in \cite{wilkinson1}). 
It was also used to analyze the large dimensional statistics of Haar matrices  in \cite{bourgade-yau}. 
Eq. \eqref{eq-overlaps} is the main tool used in the forthcoming sections to analyze the asymptotics of the overlaps in the large $N$-limit. 
Let us re emphasize the fact that the evolution equation \eqref{eq-overlaps} for $u_{i|j}(t)$
depends only on the projections of the $j$-th eigenvector $|\psi_j^0\rangle$ 
on the perturbed eigenstates $|\psi_i^t\rangle$
 and does not involve any other non-perturbed eigenvector $|\psi_\ell^0\rangle, \ell\neq j$. 
 This is a very convenient fact as we can fix a given non-perturbed eigenstate $|\psi_j^0\rangle$ and 
 work out the system of closed equations \eqref{eq-overlaps} satisfied by its $N$ 
 projections on the perturbed eigenvectors 
 $|\psi_i^t\rangle, i=1,\cdots,N$.

\subsection{Spectral density and spikes trajectories in the large $N$ limit}\label{density-trajectory}

In this section, we describe the evolution of the limiting eigenvalues density when $N\to \infty$.  
We consider the empirical spectral density of the matrix $X_t$ at time $t$ defined as 
\begin{align*}
\mu_t^N (d\lambda) := \frac{1}{N} \sum_{i=1}^N \delta_{\lambda_i(t)} (d\lambda)\,.  
\end{align*}
A classical method permits us to obtain the evolution equation for the empirical density $\mu_t^N(d\lambda)$. 
The method simply consists in computing
the infinitesimal increments over time of functional of the form
\begin{align*}
\int_\R f(\lambda) \mu_t^N(d\lambda) = \frac{1}{N} \sum_{i=1}^N f(\lambda_i(t))
\end{align*} 
where $f:\R\to \R$ is a smooth test function. This is done thanks to It\^o's formula (see for instance 
\cite[Subsection 4.3 page 248]{agz}, \cite{rogers} or more recently \cite{jp-alice} in a slightly wider context).  
With $f(\lambda)=\frac{1}{z-\lambda}$, one obtains the following Burgers evolution equation for the Stieltjes transform 
$G_N(z,t) := \frac{1}{N} \sum_{i=1}^N \frac{1}{z-\lambda_i(t)}$, 
\begin{align}\label{burgers-N}
\partial _t G_N(z,t) = - G_N(z,t) \partial_z G_N(z,t) + \sqrt{\frac{2}{\beta N}} \sum_{i=1}^N 
\frac{1}{(z-\lambda_i)^2} \frac{dB_i}{dt} + \frac{1}{2N}(\frac{2}{\beta}- 1)\partial_z^2 G_N(z,t)\,. 
\end{align} 
In the large $N$ limit, this evolution equation becomes deterministic 
%Several formulations are possible. 
and the solution is the Stieltjes transform $G$  of the limiting 
eigenvalues density $\rho(\cdot,t)$  of the matrix $X_t$ such that $\mu_t^N(d\lambda)\to \rho(\lambda,t) d\lambda$ when $N\to \infty$. 
The Stieltjes transform of the density $\rho(\lambda,t)$ is  
defined for $z\in \C\setminus\R$ as 
\begin{align*}
G(z,t) = \int_\R \frac{\rho(\lambda,t)}{z-\lambda} d\lambda\,. 
\end{align*}
This analytic function 
characterizes the probability density $\rho(\cdot,t)$ that one can compute from the imaginary part of $G$ 
near the real axis thanks
to the {\it Stieltjes inversion formula}  $\Im G(\lambda+i \varepsilon,t) \to_{\varepsilon\to 0}- \pi \rho(\lambda,t)$.  

Sending $N\to \infty$ in \eqref{burgers-N},  we see that
the dynamics of the Stieltjes transform $G$ are governed by the following {\it Burgers} evolution equation 
\begin{align}\label{burgers}
\partial_t G(z,t) = - G(z,t) \partial _z G(z,t), \quad \mbox{ with }  \quad G(z,0) = \int_\R \frac{\rho_A(\lambda)}{z-\lambda} d\lambda\,. 
\end{align}
 Interestingly, the solution of \eqref{burgers} is known \cite{shlyakhtenko} to satisfy the fixed point equation (see  also
 \cite[Proposition 4.1]{vectors2})
 \begin{align}\label{fixed-point-G}
G(z,t) = \int_0 ^1 \frac{dx}{z-a(x)- t G(z,t)}
\end{align}
where $a:[0,1]\to \R$ is the continuous function introduced in \eqref{definition-a} mapping the index $x\in [0,1]$ to the eigenvalue $a(x)$ of the matrix $A$ 
in the continuous limit $N\to \infty$.  
In the special case $A=0$, $a(x)=0$ for any $x\in[0,1]$ and the solution $G$ is fully explicit corresponding 
to the Wigner semi-circle density $\rho(\lambda,t)= \frac{1}{2\pi t} \sqrt{4t-\lambda^2}$ with radius $2\sqrt{t}$. 

One can also write the evolution equation directly in terms of the density $\rho(\lambda,t)$ itself 
by projecting the Burgers equation \eqref{burgers} on the real line thanks to the Stieltjes inversion formula: 
for $\lambda\in \R$ and $t\ge 0$,
\begin{align}\label{evolution-rho}
\partial_t &\rho(\lambda,t) + \partial_\lambda\left( v(\lambda,t) \rho(\lambda,t) \right) =0 \quad \mbox{ where } \quad v(\lambda,t) = P.V. \int_\R\frac{\rho(\lambda',t)}{\lambda-\lambda'}d\lambda' \\
  &\mbox{ and with the initial condition } \quad \rho(\lambda,0)= \rho_A(\lambda)\notag \,. 
\end{align}

Similarly as before, we will denote by $\lambda(x,t)$ the quantile function associated to the probability density 
$\rho(\lambda,t)$ such that for any $x\in [0,1]$, 
\begin{align}\label{quantile-vp}
x = \int\limits_{\lambda(x,t)}^{+\infty} \rho(\lambda,t)\, d\lambda\,. 
\end{align}
Note that if $i:=(i_N)_{N\in \N}$ is a sequence such that $i_N/N\to x\in(0,1)$, then the $i$-th eigenvalue $\lambda_i(t):=(\lambda_{i_N}(t))$ converges (almost surely) towards $\lambda(x,t)$ when $N\to \infty$. 
From \eqref{quantile-vp}, it is straightforward to check that $\partial_x \lambda(x,t)= -1/\rho(\lambda(x,t),t)$ and (using \eqref{evolution-rho}) 
that 
\begin{align}\label{time-deriv-lambda}
 \partial_t \lambda(x,t)= v(\lambda(x,t),t)\,. 
 \end{align}
We have a clear physical interpretation for the function $v(\lambda(x,t),t)$ as the speed of the particles in the scaling limit.  

The spikes trajectories become also deterministic in the large $N$ limit. We can compute them 
by  sending $N\to +\infty$ directly
in the Dyson Brownian motion equation \eqref{dyson-bm} for $j=1,\dots,\ell$. 
The limiting path of the spike $\lambda_j(t)$ for $j=1,\dots,\ell$ is
driven by the density $\rho(\cdot,t)$ satisfying \eqref{evolution-rho}. 
For $j=1,\dots,\ell$, 
\begin{align}\label{trajectory-spikes}
\dot{\lambda}_j(t) = \int_\R \frac{\rho(\lambda,t)}{\lambda_j(t)-\lambda} d\lambda \quad \mbox{ with } \quad \lambda_j(0)=a_i\,. 
\end{align}
Notice that we use the same notation for the spike trajectories $\lambda_j(t)$ for both the limiting case $N\to \infty$ and the finite dimensional case 
$N<\infty$ \footnote{In order to avoid heavy notations, we omit to use an additional super script $N$ for the eigenvalues $\lambda_i^N(t)$ at finite $N$.}. 
At the initial time $t=0$, the spike $\lambda_j(t)$ starts from a position $a_j$ outside the bulk of the spectrum 
of $A$. As $t$ increases, the spike $\lambda_j(t)$ is pushed away 
with an electrostatic force exerted by the other particles. Each particle inside the bulk of the spectrum
exerts a force which is proportional to the inverse of its distance to the spike. In such a way, the spike remains 
at a non-negative distance to the bulk at any time $t\ge 0$.

As illustrated in the next subsection, we shall nevertheless notice that the spike may be eventually 
caught back by the continuous part of the spectrum. 

%Let us now illustrate the results 
%of this section by analyzing explicitly a special and natural case.  and compute the time at which the 
%largest isolated eigenvalue gets inside the support of the continuous part).  
%
\subsection{Factor model} \label{factor-model}
In this subsection, we illustrate the results of the previous subsection by analyzing 
explicitly the special case where the matrix $A$ is of low rank $\ell$ compared to the dimension, $\ell\ll N$.
Such factor models are used in applications in biology to study population dynamics \cite{wishart} or in finance  
where the setting is nevertheless slightly different, 
see the discussion in section \ref{covariance} or \cite{review-jpb}.  
In most applications, 
the rank $\ell$ 
is fixed independently of $N$. The matrix $A$ has $\ell$ spikes $\lambda_1 \ge \lambda_2 \ge \dots \ge \lambda_\ell$ and 
$0$ is an eigenvalue of $A$ with multiplicity $N-\ell\sim N$. 
The structure of the matrix $A$ is therefore very simple with only a few relevant factors that one wants to estimate.  

The few spikes do not bring any macroscopic contribution to the empirical density $\rho(\cdot,t)$ of the particles
and in the large $N$ limit, we recover the Wigner semicircle density centered at $0$ with radius $2\sqrt{t}$, 
\begin{align}\label{semicircle-t}
\rho(\lambda,t) = \frac{1}{2\pi t} \sqrt{4t-\lambda^2}\,, \quad -2\sqrt{t} \le \lambda \le 2\sqrt{t}\,. 
\end{align}   
The speed of the particles inside the spectrum can be computed explicitly as well: it is linear given for $t>0, |\lambda| \le 2\sqrt{t}$ by
\begin{align*}
v(\lambda,t) = \frac{\lambda}{2 t} \,. 
\end{align*}

\begin{figure}\center
\includegraphics[scale=0.9]{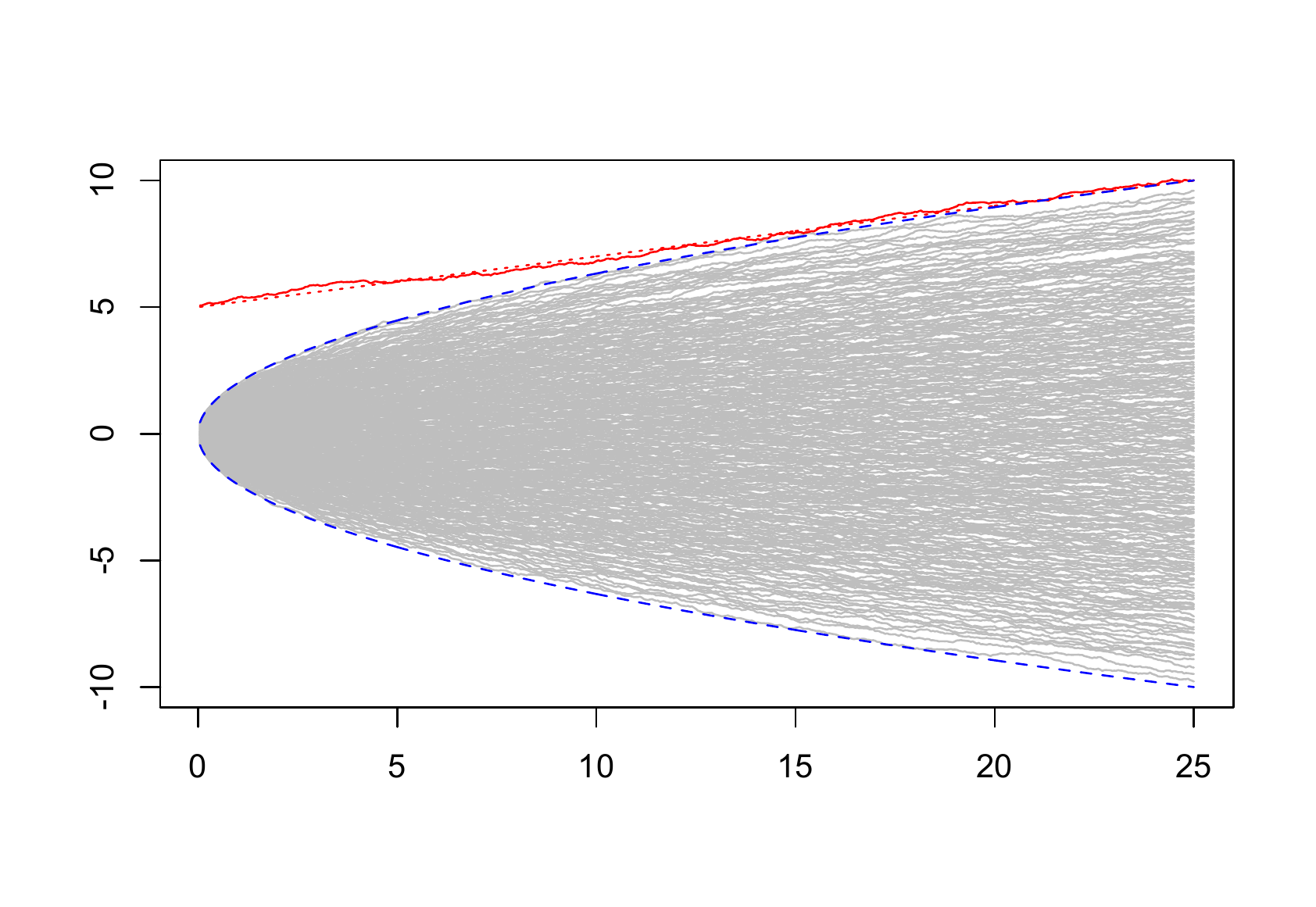}
\caption{(Color online). Sample trajectories of the eigenvalues of the matrix $X_t$ defined in \eqref{def-X} 
where $A$ has only one non zero eigenvalue $a_1=5$, as a function of time $t>0$. 
The grey lines represent the 
eigenvalues trajectories in the bulk. The blue dashed lines 
represent the trajectories of the edges $\pm2\sqrt{t}$. 
The red plain line is the sample path of the spike $\lambda_1(t)$ and the red dashed line is $5+t/5$.
Beyond $t_c=25$, the spike eigenvalue $\lambda_1(t)$  is 
``swallowed'' by the Wigner sea and disappears.}\label{fig-eigenvalue-factor-model}
\end{figure}

With such a simple form \eqref{semicircle-t} for the limiting density of particles, it turns out that 
the ordinary differential equation \eqref{spike1} can be solved explicitly thanks to elementary computations. 
We obtain for any $j=1,\dots,\ell$, 
\begin{align*}
\lambda_j(t) = a_j + \frac{t}{a_j}\,. 
\end{align*}
Comparing this value of the $j$-th spike with the value of the edges of the spectrum at time $t$, 
we easily check that for any $a_j\neq 0$, the bulk eventually catches up the isolated particle 
$\lambda_j(t)$ at the critical time $t_c^j=a_j^2 $ at which $\lambda_j(t_c^j)=2\sqrt{t_c^j}$, beyond which the spike is 
``swallowed'' by the Wigner sea and disappears. See Fig. \ref{fig-eigenvalue-factor-model} for an illustration 
of a sample path of the eigenvalues of $X_t$ when the initial matrix $A$ has rank one.

\section{Eigenvector in the bulk of the spectrum}\label{bulk}

In the bulk of the spectrum, the mean spacings $\delta a$ 
between the eigenvalues of the matrix $A$ is approximately of order $1/N$ 
and depends on the position $a$ in the spectrum and the local density $\rho_A(a)$ of particles near $a$ as
$ \delta a \sim 1/(N\rho_A(a))$. 

When one perturbs the initial matrix $A$ by adding the matrix $H_t$, it is well known that 
one should compare the magnitude of the entries of the perturbation $H_t(ij) \sim \sqrt{t}/\sqrt{N}$
with the mean level spacing $\delta a$ of the non-perturbed matrix $A$. 

There are therefore three distinct regimes of perturbation which lead to different asymptotics for the relation 
between the perturbed and non-perturbed eigenstates: 
 \begin{itemize}
 \item The microscopic or perturbative regime corresponds to values of $t:=t_N$ depending on $N$ such that 
 \begin{align*}
t_N\ll \frac{1}{N}\,. 
\end{align*}
For such values of $t:=t_N$, the perturbation matrix $H_t$ is in fact asymptotically small compared to $A$ and for any fixed $i$, 
the 
eigenvector $|\psi_j^t\rangle$ of $X_t$ converge to those of $A$ when $N\to \infty$ in the $L^2$ norm,
\begin{align}\label{conv-psi_i}
||\psi_j^{t_N} - \psi_j^0||_2 \longrightarrow 0\,. 
\end{align}
One can even obtain an asymptotic expansion for $|\psi_j ^t\rangle$ around $|\psi_j ^0\rangle$ using \eqref{evolution-psi_i}.
This regime is rather trivial and will not be further considered  in this paper. 
 \item The second mesoscopic 
 regime establishes a smooth crossover between the microscopic and macroscopic regimes. 
  It corresponds to values of  $t:=t_N$ which are inversely proportional to $N$ i.e. such that 
  there exists $\tau >0$ fixed such that
\begin{align*}
t_N = \frac{\tau}{N\rho_A(a_j)}\,. 
\end{align*}
Although the operator norm 
of the matrix $H_{t_N}$ tends to $0$ when $N\to \infty$, this regime is non-perturbative in the sense that we do not have the convergence 
\eqref{conv-psi_i} of $|\psi_i^t\rangle$ towards $|\psi_i ^0\rangle$. 
This non trivial rotation of the eigenvectors may appear surprising at first sight 
(it is generated by the addition of a microscopic perturbation) but
 is in fact simply due to the small spacings $\delta a$ between the eigenvalues of $A$ in the bulk of the spectrum.  
We shall analyze this regime in details in subsection \ref{crossover} thanks to the overlaps equation \eqref{eq-overlaps}.
We will see that, for any $i,j$ such that $\lambda_i,a_j$ lie in the bulk of the spectrum, the vector $|\psi_j^0\rangle$
is {\it localized} in a finite dimensional subspace of $\R^N$ (the dimension is proportional to $\tau$ in the limit $1\ll \tau\ll N$, see below) 
spanned by the perturbed eigenvectors 
$|\psi_i^t\rangle, j \sim i$ associated to the finitely many closest neighboring eigenvalues 
of $a_j$ such $N (a_j -\lambda_i(t) ) \sim 1$, as $N\to \infty$.  
 \item The macroscopic regime corresponds to values of $t=O(1)$ which do not depend on $N$.
Even though the perturbation is macroscopic, we are still able to extract  information 
on the non-perturbed eigenstate $|\psi_j ^0\rangle=|\phi_j\rangle$
from the observation of the perturbed eigenstates $|\psi_i^t\rangle$ for general $i,j$. 
Indeed we compute explicitly the asymptotic mean overlaps $\left[  \langle\psi_i^t|\psi_j^0\rangle^2\right]$ 
which are proportional to $1/N$ in the large $N$ limit, 
using again the overlap equation \eqref{eq-overlaps}. If 
$i:=(i_N)_{n\in \N}$ and $j:=(i_N)_{n\in \N}$ are sequences such that $i_N/N \to x \in(0,1)$ and 
$j_N/N\to y\in(0,1)$ when $N\to \infty$,  our result reads (see below)
\begin{align*}
\left[ \langle\psi_i^t|\psi_j^0\rangle^2\right] \underset{N\to \infty}{\sim} \frac{1}{N} \, u(x,y,t) 
\end{align*}
where the function $u(x,y,t) \sim O(1)$ is determined explicitly for any matrix $A$ in terms of the  trajectory of the limiting 
density $(\rho(\cdot,s))_{0\le s\le t}$  
described in \eqref{evolution-rho}. The function $u(x,y,t)$ can be computed explicitly  in terms of elementary 
functions in some special cases. 
 \end{itemize}

\subsection{Crossover regime} \label{crossover}
In this section, we consider the crossover regime $t\propto\tau/N$ where $\tau>0$ is fixed. 
\subsubsection{General limiting random equation}
Let us first fix a sequence of indices $j:=(j_N)_{N\in \N}$ such that the eigenvalue $a_j:=a_{j_N}$ converges when $N\to \infty$
to a fixed position $a(y)$ in the bulk of the spectrum. The sequence $(j_N)_{N\in \N}$ satisfies the asymptotic scaling relation
\begin{align}\label{scaling-j}
 \frac{j_N}{N} \underset{N\to \infty}{\longrightarrow} y.
\end{align}
In the following, we omit to write the dependence in $N$ for the subscript $j$. 

Recall the overlap equation on the projections of the state $|\psi_j^0\rangle$ on the perturbed eigenstates $|\psi_i^t\rangle$ at time $t$ 
may be written
  \begin{align}\label{overlap-eq-2}
\partial_t\, u_{i|j} (t) = \frac{1}{N} \sum_{k \neq i} \frac{u_{k|j}( t)- u_{i|j}( t)}{(\lambda_k(t)-\lambda_i(t))^2}\quad 
\mbox{ with } \quad u_{i|j}(0) = \delta_{ij}\,.
\end{align} 

We first need to establish the convergence of the eigenvalues process $(\lambda_k(t))_{t\ge 0}$ satisfying \eqref{dyson-bm} 
on the microscopic time scale $t=\tau/(N\rho_A(a_j))$. 
From the Dyson equation \eqref{dyson-bm}, we see that the re-scaled eigenvalues process 
converges weakly, when $N\to \infty$, towards an {\it infinite} random particle system $(x_k(\tau))_{k \in \Z,\tau \ge 0}$  such that 
\begin{align}\label{eq-x_k}
dx_k(\tau) = \sqrt{\frac{2\rho_j}{\beta}} \, dB_k(\tau) + \rho_j \sum_{\ell \in \Z\setminus\{ k\} } \frac{d\tau}{x_k(\tau)-x_\ell(\tau)}  
\quad \mbox{ with } \quad  x_k(0) =  k  
\end{align}
where $\rho_j\equiv \rho_A(a_j)>0$. 
Such diffusion processes in the space of real sequences indexed by $\Z$ 
 have been very recently defined \cite{osada} through a limiting procedure where the infinite sum over $\Z$ is truncated.  
The process $(x_k(\tau))_{k \in \Z,\tau \ge 0}$ is invariant by translation at any time and stationary in time. 
 
The convergence of the process $(\lambda_k(\frac{\cdot}{N\rho_j}) )_{1\le k\le N}$ towards the process $(x_k(\cdot))_{k \in \Z}$ reads
\begin{align}\label{conv-sine-dyn}
\left(N \rho_j\, \left ( \lambda_{j+k}(\frac{\tau}{N\rho_j})  - a(y) \right)\right)_{ 1-j \le k \le N-j, \tau \ge 0} 
\longrightarrow \quad (x_k(\tau))_{k \in \Z,\tau \ge 0}\,,
\end{align} 
where $a(y)$ is the limiting position of the eigenvalue $a_j$ at the initial time. 
Note that we have to shift the indices by $j$ lags before sending $N$ to $\infty$. 

The convergence \eqref{conv-sine-dyn} describes the joint local statistics in law of the eigenvalues process 
in the microscopic 
vicinity (i.e. in a region of width $1/N$) of the $j$-th eigenvalue when one re-scales time by a factor $1/N$. This convergence is usually written at a fixed time where one obtains
the classical limiting Sine$_\beta$ point process (for $\beta=2$, it is the Sine kernel determinantal point process; see 
also \cite{valko-sine} or \cite{laure-sine}
for a description of the law of the Sine$_\beta$ point process for general $\beta>0$). 
The convergence \eqref{conv-sine-dyn} is valid for a sequence of indices $j$ proportional to $N$
as in \eqref{scaling-j} where $y$ is the limiting position of the $j$-th eigenvalue at the initial time. 
There are no constraints on the index $k$ and in the scaling limit the process $(x_k)$ is indexed on $\Z$. Note that the eigenvalues 
standing at a macroscopic (or mesoscopic) distance away from the $j$-th eigenvalue at the initial time are ejected in $\pm \infty$ 
in the scaling limit.

The convergence \eqref{conv-sine-dyn} at the initial time $\tau=0$ 
is obtained thanks to the smooth allocation of the eigenvalues $a_k=a(\frac{k}{N})$ of $A$ 
such that $N\rho(a_j) (a(\frac{k}{N}) - a(\frac{j}{N})) \sim  k - j$. 
The re-scaling in space by the factor $N \rho_A(a_j)$ is chosen to have a level spacing 
 approximately equal to $1$ at the initial time. At later time $\tau >0$, the mean level spacing in between the $x_k(\tau)$ 
 remain macroscopic  because of
 the electrostatic repulsion between the particles in the system (which is sufficiently strong to prevent any collision \cite{agz}).  

Now we want to use the convergence of the eigenvalue process to study the overlap equation
\eqref{overlap-eq-2} in the double scaling limit $t=\tau/N, N\to \infty$. 

Recall that the sequence $j$ is fixed and let us introduce the finite family of rescaled overlaps (we drop heretoforth the explicit 
dependence on $j$ which will only appear through a rescaling of time):
\begin{align*}
v_i( \tau) =   u_{j+i|j}(\frac{\tau}{N\rho_j}), \quad i = 1-j ,\dots,N-j.
\end{align*} 
We denote by $\ell^1 (\Z)$ the space of real sequences $(v_i)_{i \in \Z}$ indexed by $\Z$ such that  $\sum_{i \in \Z} v_i < +\infty$. 

Using \eqref{overlap-eq-2} and \eqref{conv-sine-dyn}, we can now determine the limiting evolution equation governing the dynamics of 
the infinite sequence $(v_i(\tau))_{i \in \Z}$ in the space $\ell^1(\Z)$   
as a function of the time $\tau$ in the scaling limit $N\to \infty$. 
This evolution equation obtained from Eq. \eqref{overlap-eq-2} may be written in terms of the diffusion process $(x_k(\tau))_{k \in \Z,\tau \ge 0}$,
for $i \in \Z$, as
\begin{align}\label{crossover-general-eq}
\partial_\tau \, v_i(\tau) = \rho_j \sum_{k \in \Z \setminus \{i\}} \frac{v_k(\tau)- v_i(\tau)}{(x_{k}(\tau)-x_{i}(\tau))^2}\quad \mbox{ with } 
\quad v_i(0) = \delta_{i0} \, .
\end{align}
Recall that the process $(x_k(\tau))_{k \in \Z,\tau\ge 0}$ (satisfying \eqref{eq-x_k}) is stochastic so that 
the process $(v_i(\tau))_{i \in \Z,\tau \ge 0}$ of the overlaps is also stochastic (even though we took a limit $N\to \infty$). 
 
Equation \eqref{crossover-general-eq} describes the evolution in time of the sequence $v_i(\tau)$  
in the double scaling limit $ t=\tau/N, N\to \infty$. Note that the limit point $v(\tau) := (v_i(\tau)) _{i \in \Z}$ satisfies  
$\partial _\tau \sum_{i \in \Z} v_i(\tau)= 0$ and hence $  \sum_{i \in \Z} v_i(\tau)= 1$ for any $\tau\ge 0$ as it should be. 
The 
sequence $(v_i(\tau))_{i \in \Z}$ can also be seen as a probability distribution on $\Z$. 
At the initial time $\tau=0$, this probability distribution 
is a Dirac delta function (the perturbation is null and we measure the eigenvectors perfectly). As 
$\tau$ grows, this {\it random} distribution on $\Z$ broadens according to the dynamical equation \eqref{crossover-general-eq}.   

As a conclusion of this paragraph, we conjecture (a few arguments used here are heuristical) 
 the following convergence in law in the space of continuous process in $\R^n$, 
for any $n\in \N$, $T>0$, $j$ scaling with $N$ according to \eqref{scaling-j} and any $i_1,i_2,\dots, i_n \in \Z$,  
\begin{align}
\left(\left[ \langle \psi_{j + i_1}^{\tau/(N\rho_j)}  |\psi_j^0\rangle^2\right] , \left[ \langle \psi_{j + i_2}^{\tau/(N\rho_j)}  |\psi_j^0\rangle^2\right], \dots,
\left[ \langle \psi_{j + i_n}^{\tau/(N\rho_j)}  |\psi_j^0\rangle^2\right] \right)_{0 \le \tau \le T }\notag \\
\Rightarrow 
\left( v_{i_1}(\tau), v_{i_2}(\tau), \dots, v_{i_n}(\tau)\right)_{0\le\tau \le T}\label{conjecture-conv1}
\end{align}
where the stochastic process $(v_i(\tau))_{i \in \Z,\tau \ge 0} $ taking values in $\ell^1(\Z)$  satisfies the limiting equation \eqref{crossover-general-eq}.
This convergence along the finite dimensional marginals can be translated in a convergence in law in the space of continuous process 
in $\ell^1(\Z)$ of the sequence $([  \langle \psi_{j + i}^{\tau/(N\rho_j)}  |\psi_j^0\rangle^2])_{i \in \Z, 0 \le t \le T}$ 
towards $(v_i(\tau))_{i \in \Z, 0 \le t \le T}$ (the tightness is obvious as we work with bounded sequences in $\ell^1(\Z)$: the sum of the elements is $1$). 
%Another interesting question would be to prove the weak convergence of the local density of state $j$, i.e. of the random measure process
%\begin{align*}
%\left(\sum_{i \neq j} \E_{\lambda}[ \langle \psi_{j + i}^{\tau/N}  |\psi_j^0\rangle^2] \, \delta_{\lambda_{j+i}(\frac{\tau}{N})}(dx)\right)_{0 \le \tau \le T}
%\end{align*}
%towards $\left(\sum_{i \in \Z} v_i(\tau) \delta_{x_i(\tau)}(dx) \right)_{\tau \ge 0}$. 

It would also be  interesting to obtain precise insights on the fluctuations of the stochastic process 
$(v_i(\tau))_{i\in \Z,\tau \ge 0}$. In particular, the heat kernel associated to the stochastic equation
\eqref{crossover-general-eq} is very intriguing. A very nice result would be to compute the mean heat kernel
in the stationary case where the initial distribution of the point process $(x_k(0))_{k \in \Z}$ is the 
Sine$_\beta$ law \cite{question-bourgade} (see \cite{valko-sine} or \cite{laure-sine} for a reminder on the Sine$_\beta$ point processes).   

In the following, we compute explicitly the deterministic sequence $(v_i(\tau))_{i\in \Z,\tau \ge 0}$ associated to the non stochastic 
but most probable 
trajectory of the particle system $(x_i(\tau))_{i \in \Z, \tau \ge 0}$. 

\subsubsection{The most probable deterministic evolution associated to Fekete trajectories.}

The most probable path for the infinite dimensional diffusion process $(x_k(\tau))_{k\in \Z}$, also called the {\it Fekete} trajectory,
satisfies \eqref{eq-x_k} where the noise terms $dB_k, k \in \Z$ have been set to $0$. 
One can check that the deterministic Fekete trajectories of the $x_k$ are in fact constant in time such that, for all $\tau \ge 0$,
\begin{align*}
x_k(\tau) = k \,. 
\end{align*} 
With those Fekete trajectories, the equation  on the overlaps \eqref{crossover-general-eq} becomes deterministic 
\begin{align}\label{crossover-deterministic-eq}
\partial_\tau \, v_i(\tau) = \rho_j \sum_{k \in \Z \setminus \{i\}} \frac{v_k(\tau)- v_i(\tau)}{(k-i )^2}\quad \mbox{ with } 
\quad v_i(0) = \delta_{i0} \,.
\end{align}
Setting $v(\tau) := (v_i(\tau)) _{i \in \Z}$,  Eq. \eqref{crossover-deterministic-eq} can be rewritten as 
\begin{align}\label{eq-with-U}
\partial_\tau v = - \mathcal{U} v(\tau)
\end{align}
where $\mathcal{U}$ denotes the linear operator in $\ell^1(\ZZ)$ such that for $v:=(v_i)_{i \in \Z} \in \ell^1(\Z)$,  
\begin{align}\label{def-op-U}
(\mathcal{U}v)_i:=  \rho_j \sum_{k \neq i} \frac{v_i-v_k}{(i-k)^2} \,. 
\end{align}
Eq. \eqref{eq-with-U} was also solved en passant in \cite{bourgade-univ}. 
We revisit the derivation proposed in \cite{bourgade-univ}. 

Following \cite{bourgade-univ}, we introduce the Fourier transform defined for a sequence $(v_i)_{i \in \ZZ}\in \ell^1(\ZZ)$ and
$\xi\in \R$ as 
\begin{align*}
{\hat v}(\xi) := \sum_{k \in \ZZ} e^{-i 2\pi\xi k} v_k\,. 
\end{align*}
We can easily compute the Fourier transform of the sequence $\mathcal{U}v$ as a function of $\hat{v}$
\begin{align*}
\widehat{\mathcal{U}v}(\xi) &= \rho_j  \sum_{k \in \ZZ}  e^{-i 2\pi\xi k} \sum_{j\neq k} \frac{v_k-v_j}{(k-j)^2} \\
&=  \rho_j  \sum_{\ell\neq 0} \frac{1}{\ell^2} \sum_{k \in \ZZ}e^{-i 2\pi\xi k} (v_k - v_{k-\ell}) \\
&= \rho_j \hat{v}(\xi) f(\xi)
\end{align*}
where 
\begin{align*}
f(\xi)= \sum_{\ell\neq 0} \frac{1}{\ell^2}  (1- e^{-i 2 \pi \xi \ell})\,. 
\end{align*}
For $\xi\in [0,1]$, we have the explicit form 
\begin{align*}
f(\xi)= 2\pi^2 \xi (1-\xi)\,. 
\end{align*}
This explicit form for the function $f$ is different  from the one proposed in \cite{bourgade-univ} where the quadratic term
appears to be missing.  
We can now solve \eqref{eq-with-U} in Fourier space. For $\xi\in [0,1]$,  
\begin{align*}
\hat{v}(\xi,\tau) = \exp\left(-2 \pi^2 \xi (1-\xi)  \rho_j  \tau \right) 
\end{align*}
 where we have used the initial condition $v_i(0)=\delta_{i0} $ (i.e. $\hat{v}(\xi,0)=1$). 
 
Now we recover the sequence $v_n(\tau)$ through the inverse Fourier transform 
\begin{align}\label{solution-v}
v_n(\tau) = \int_0^1 \exp\left(-2 \pi^2 \xi (1-\xi)  \rho_j  \tau \right) \cos (2\pi\xi n)  d\xi\,. 
\end{align} 

 \subsubsection{Localization of the non-perturbed eigenstate $|\psi_j^0\rangle$ in a cone of dimension $\propto\tau$} 
In order to estimate the dimension of the subspace containing $|\psi_j^0\rangle$ (associated to the non-perturbed eigenvalue $a_j$)
generated by the eigenvectors $|\psi_i^{\tau/(N\rho_j)}\rangle$ associated to 
the neighboring eigenvalues $\lambda_i(\frac{\tau}{N})\sim a_j$ for $i \sim j$, we study
the asymptotic of the overlap $v_n(\tau)=\lim_{N\to \infty} \langle \psi_{j+n}^{\tau/(N\rho_j)}|  \psi_j^0\rangle ^2$ 
in the double scaling limit $\tau\to +\infty$ with $n$ scaling with $\tau$ as $n = p\tau$ 
where $p>0$ is a fixed parameter. In this regime, one has:
\begin{align}
v_n(\tau)& =2 \int_0^{1/2}  \exp\left(- 2 \pi^2 \tau  \xi (1-\xi) \rho_j \right)   \cos(2\pi\xi n) d\xi\notag \\
&= \frac{2}{\tau} \int_0^{\tau/2} \exp\left(- 2 \pi^2 x(1-\frac{x}{\tau}) \rho_j \right)\cos(2\pi x\frac{n}{\tau} ) dx \notag \\
&\sim \frac{2}{\tau} \int_0^{+\infty} \exp(-2 \pi^2 \rho_j x) \cos(2\pi x p ) dx = \frac{1}{\tau} \frac{\rho_j}{p^2+\pi^2 \rho_j^2}\notag\\& =  \frac{\tau \rho_j}{ n^2 +\pi^2 \tau^2 \rho_j^2 } = \frac{1}{\rho_j \tau} F \left(\frac{n}{\rho_j\tau}\right)\label{cauchy}
\end{align}
where $F(u)=1/( u^2+\pi^2)$ is the (normalized) Cauchy distribution.
The equivalence in the third line holds in the double scaling limit $\tau,n\to \infty$ with 
 $n/\tau=p$, $p>0$ fixed. 
This computation proves that, if $\tau>0$ is large, the state $|\psi_j^0\rangle$ is (almost) entirely contained in the subspace spanned by the 
eigenvectors $|\psi_{j + n}^{\tau/(N\rho_j)}\rangle$ with indices $j+n$ 
such that $n$ is proportional (or smaller) to $\tau$. The eigenvectors $|\psi_{j + n}^{\tau/(N\rho_j)}\rangle$ such that $n$ is much larger than $\tau$ have 
very small overlaps with $|\psi_j^0\rangle$. 
 
The above ``Cauchy-flight'' shape for the diffusion of the overlaps in the mesoscopic regime 
$\tau \gg 1$ (i.e. $1/N\ll t\ll 1$) was already found in \cite{wilkinson1} -- see also \cite{vectors2}. 
The authors of \cite{wilkinson1} solve the general overlap equation \eqref{eq-overlaps} (see Eq. (4.7) in their paper) in the mesoscopic 
regime through a non rigorous computation. They use the rigidity of the eigenvalues which basically reduces to working with the Fekete trajectories 
as we do in this subsection, but they also need to approximate the overlap equation \eqref{eq-overlaps}  
with a continuous equation where the sum 
is replaced by an integral (this approximation is not true in the mesoscopic regime but only 
in the macroscopic regime $t\sim 1$, see the next section \ref{third-regime}). They finally take an extra assumption, 
which is physically sound but not justified, 
on the kernel $R$ (see Eq. (4.9) in \cite{wilkinson1}) involved in their continuous equation (Eq. (4.8) of \cite{wilkinson1}). 

Our approach is (to our eyes) more transparent as we solve the overlap equation \eqref{crossover-deterministic-eq} in the microscopic regime 
 $t=\tau/(N\rho_j)$ when the eigenvalues follow the Fekete trajectories. Our solution is explicit given in Eq. \eqref{solution-v}
 and we finally 
re obtain the Cauchy shape \eqref{cauchy} for the local density of the state $|\psi_j^0\rangle$ going
from the microscopic to  
the mesoscopic regime by sending the parameter $\tau$ to $+\infty$.

\subsection{Non-perturbative regime} \label{third-regime}
We now consider the case where $t>0$ is fixed independently of $N$. In this regime, we expect the distribution 
of the overlaps to be much more spread out compared to the other regimes: the non-perturbed eigenstates 
are delocalized in the basis of the perturbed eigenvectors. All the mean squared 
overlaps have the same order of magnitude of order $1/N$ for large $N$.   

We start by deriving the evolution equation of the local density of the state $|\psi_j^0\rangle$ for a fixed index $j$. 

\subsubsection{Local density of state}

The local density of the state $|\psi_j^0\rangle$ describes the allocation of the mean squared projections of the 
non-perturbed state $|\psi_j^0\rangle$
on the basis of the perturbed eigenvectors $|\psi_i^t\rangle$. It is a probability measure defined as 
\begin{align*}
\nu_N^{(j,t)}(d\lambda):= 
\sum_{i =1} ^N [\langle\psi_j^0|\psi_i^t\rangle^2] \, \delta_{\lambda_i(t)}(d\lambda) = \sum_{i =1} ^N u_{i |j}(t) \delta_{\lambda_i(t)}(d\lambda)\,. 
\end{align*} 
Let us denote by $U_N(z,t)$ the Stieltjes transform of this probability measure 
\begin{align*}
U_N^{(j)}(z,t) := \int_\R \frac{\nu_N^{(j,t)}(d\lambda)}{z-\lambda}= \sum_{i=1}^N \frac{u_{i|j}(t)}{z-\lambda_i(t)}\,. 
\end{align*}
It is easy to check that $U_N^{(j)}(z,t)$ is equal to the local resolvent in the sense that, for any $j$ and $z\in \C\setminus \R$, we have  
\begin{align*}
U_N^{(j)}(z,t)  =\langle \psi_j^0 |(z- X_t)^{-1} |\psi_j^0\rangle \,. 
\end{align*}
Now, using the Dyson equation for the eigenvalues \eqref{dyson-bm} and the overlap evolution equation \eqref{eq-overlaps}, 
we obtain (using again It\^o's Formula) the following evolution equation for the local resolvent 
\begin{align}\label{burgers-local}
\partial _t U_N^{(j)}(z,t) &= - G_N(z,t) \partial_z U_N^{(j)}(z,t) + \sqrt{\frac{2}{\beta N}} \sum_{i=1}^N 
\frac{u_{i|j}(t)}{(z-\lambda_i)^2} \frac{dB_i}{dt} + \frac{1}{2N}(\frac{2}{\beta}- 1)\partial_z^2 U_N(z,t), \\
U_N^{(j)}(z,0)& = \frac{1}{z-a_j} \notag
\end{align}
where $G_N(z,t):= \frac{1}{N} \sum_{i=1}^N \frac{1}{z-\lambda_i(t)}$ satisfies the Burgers equation \eqref{burgers-N}.  
We give the proof of Formula \eqref{burgers-local} in  Appendix \ref{burgers-local-proof}.  
 Summing over $j=1,\dots,N$ in \eqref{burgers-local}, one actually recovers the Burgers equation \eqref{burgers-N}
 describing the evolution of the Stieltjes transform $G_N$. 
 
It is not difficult to see that the stochastic partial differential equation \eqref{burgers-local} also becomes deterministic 
in the large $N$-limit. We denote by $x\in (0,1)$ the limit point of the sequence $(j_N/N)_{N\in \N}$, by $a(x)$ the limiting 
position of the $j$-th eigenvalue at time $0$ and  
by $U(z,a(x),t)$ the limiting value of $U_N^{(j)}(z,t)$ when $N\to +\infty$. 
The equation on the limiting local resolvent $U(z,a(x),t)$ reads 
\begin{align}\label{limiting-local-resolvent-evolution}
\partial_t U(z,a(x),t) = - G(z,t) \partial_z U(z,a(x),t),  \quad \mbox{ with } \quad U(z,a(x),0) = \frac{1}{z-a(x)}
\end{align} 
where $G$ satisfies the limiting Burgers equation \eqref{burgers}. 
 Recalling the fixed point equation Eq. \eqref{fixed-point-G} satisfied by $G(z,t)$, 
 it is easy to check that the solution of $U(z,a(x),t)$ such that 
 \begin{align*}
G(z,t) = \int_0^1 U(z,a(x),t) dx
\end{align*}
is actually given, for any $x\in (0,1),z\in \C\setminus \R, t\ge 0$ by 
\begin{align}\label{solution-U}
U(z,a(x),t) = \frac{1}{z-a(x)-t G(z,t)}\,. 
\end{align} 
 This explicit solution of \eqref{limiting-local-resolvent-evolution} is quite remarkable. 
 
 We have established the (almost sure) weak convergence 
 of the local density of the state $|\psi_{j_N}^0\rangle$ where $j_N$ is a sequence such that $j_N/N\to x$  
 towards the unique probability measure whose Stieltjes transform is given by the holomorphic function $U(\cdot,a(x),t)$
 given in \eqref{solution-U}. 
 The limiting local resolvent $U(z,a(x),t) $ given in \eqref{solution-U} was already obtained by Shlyakhtenko in \cite{shlyakhtenko} 
 using Free probability theory. We think the Dyson style approach developed here is very intuitive, shedding new lights on this result.  
% extending the Wigner semicircle Th 
% The proof of this result may be seen as the counterpart for
% the eigenvectors to Shlyakhtenko's result \cite{shlyakhtenko} 
% extending the Wigner semicircle Theorem 
% to random matrices of the form $X_t=A+H_t$. 

\subsubsection{Continuous equation for the overlaps}
As $N\to \infty$, the limiting overlaps are described in terms of the continuous function $u(\cdot,y,t):[0,1]\to\R_+$ 
such that 
\begin{align}\label{relation-u-U}
\int_0^1 \frac{u(x,y,t)}{z-\lambda(x,t)} dx = U(z,a(y),t) \,. 
\end{align}
The Cauchy problem satisfied by the function $u(\cdot,y,t)$ 
can be determined sending $N\to\infty$ in the discrete Eq. \eqref{eq-overlaps}; We obtain 
 \begin{align}\label{overlap-eq-non-perturb}
\partial_t \, u(x,y,t) = P.V. \int_0^1 \frac{u(z,y,t)-u(x,y,t)}{(\lambda(z,t)-\lambda(x,t))^2}   dz, \quad u(x,y,0) = \delta(x-y)
\end{align}
 where the $\lambda(x,t)$ describes the limiting path of the eigenvalue with index $x:= \lim i_N/N$.
 
If $i:=(i_N)$ and $j:=(j_N)$ are two sequences such 
that $i_N/N\to x$ and $j_N/ N\to y$ where $x,y\in (0,1)$, 
then we have the following convergence of the overlap 
\begin{align}\label{conv-overlaps-bulk}
N u_{i|j}(t) \underset{N\to \infty}{\longrightarrow} u(x,y,t)
\end{align}
 where $u$ is the unique solution of the Cauchy problem 
\eqref{overlap-eq-non-perturb}.

In terms of the local density of the state $|\psi_j^0\rangle$ indexed this time with the indices $i=1,\dots,N$, 
the convergence Eq. \eqref{conv-overlaps-bulk} is equivalent to the almost sure weak convergence 
\begin{align*}
 \sum_{i =1}^N  \left[ \langle\psi_i^t |\psi_j^0 \rangle^2 \right] \, \delta_{i/N} (dx)  \Rightarrow u(x,y,t) \, dx   
\end{align*}
where $u$ is the solution of \eqref{overlap-eq-non-perturb} and where $y:=\lim_{N\to \infty} j_N/N$.

 It is useful to index the overlap as a function of the eigenvalues instead of the indices. We denote by $w$ the 
function such that if $\lambda:= \lambda(x,t),\mu:= \lambda(y,0) , x,y\in(0,1)$ and $t\ge 0$,  
\begin{align}\label{def-w}
w(\lambda,\mu,t) := u(x,y,t)\,. 
\end{align}  

As a consequence of Eq. \eqref{solution-U}, 
it is plain to deduce (using the relation \eqref{relation-u-U}) that the solution $u(x,y,t)=w(\lambda(x,t),\lambda(y,0),t)$ 
of the Cauchy problem \eqref{overlap-eq-non-perturb} 
 is given in terms of $\rho(x,t)$ 
 and of its Hilbert transform $v(x,t)$ (defined in \eqref{evolution-rho}) as 
 \begin{align}\label{sol-non-perturb}
u(x,y,t) =w(\lambda,\mu,t) = \frac{t}{(\lambda-t \, v(\lambda,t)-\mu)^2 + t^2 \pi^2 \rho(\lambda,t)^2},
\end{align}
where we have used the same short hand notation $\lambda:= \lambda(x,t)$ and $\mu:= \lambda(y,0)$. 

It is easy to check that the solution $u(x,y,t)$ of \eqref{overlap-eq-non-perturb} satisfies $\int_0^1 u(x,y,t) dx = 1$ for any $t\ge 0, y\in [0,1]$.
The function $w$ satisfies $\int_\R w(\lambda,\mu,t) \rho(\lambda,t) d\lambda=1$ for any $t\ge 0, \mu\in \R$.

 It would be interesting to generalize this explicit function expressed 
 in terms of the density of particles $\rho(\cdot,t)$ and the velocity field $v(\cdot,t)$ (which in the present 
 case is equal to the Hilbert transform of the probability density $\rho(\cdot,t))$) 
 in a large deviation regime where $\rho(\cdot,t),v(\cdot,t)$ obey the Euler-Matytsin equations (see for example \cite{joel}).  
 This would allow one to establish further interesting connections between the Harish-Chandra-Itzykson-Zuber integral 
 and free probabilities (see \cite{joel} and \cite{crossover-addition}).

As a challenging open problem, we think it would be interesting to further characterize the fluctuations  
of the family of overlaps $\{ \sqrt{N} \, \langle \psi_j^0 |\psi_i^t\rangle, i =1,\dots,N\}$
 for a fixed value of $j$ scaling with 
$N$ as in \eqref{scaling-j}, in the limit of large $N$. 
Denoting by $i_1<\dots<i_n\in \{1,\dots,N\}$ a finite subset of indices such that $i_k/N \to x_k$ for some 
$x_k\in \R$, $k=1,\dots,n$, 
we conjecture that the sub family of random variables
\begin{align*}
\left(\sqrt{N} \, \langle \psi_j^0 |\psi_{i_k}^t\rangle\right) _{k =1,\dots,n}
\end{align*}
converge in law to a centered Gaussian vector $(g_1,\dots, g_n)$ whose entries $g_k$ are independent with respective variance 
$\langle g_k^2 \rangle = u(x_k,y,t)$ where $j_N/N\to y$.
In section \ref{gaussian-fluctuations}, we are able to solve a related problem on the limiting Gaussian 
fluctuations of the overlap between the perturbed and non-perturbed eigenvector  associated to an 
isolated eigenvalue (standing away from the bulk density at the initial time) with a moment method.

\subsubsection{Perfect matching with the mesoscopic regime}

It is interesting to note the perfect matching between formulas \eqref{sol-non-perturb}
and \eqref{cauchy} in the mesoscopic regime. Indeed, applying Taylor formula when $t\to 0$, we have 
\begin{align*}
\lambda(x,t)-t \, v(\lambda(x,t),t)-\mu &= \lambda(x,0) + t \partial_t \lambda(x,0) - t v(\lambda(x,0),0)-\mu + O(t^2)\\
&= \lambda(x,0) - \mu + O(t^2)
\end{align*}
where we have used \eqref{time-deriv-lambda} to obtain the second line. The speed term cancels at the first order in $t$ when $t$ is small.  
 Therefore, if $t=\tau/(N\rho_j)$ and $x= i/N, y =j/N$ where $n:=i-j \in \Z$ is fixed independently of $N$, then  
 \begin{align*}
\lambda(x,\frac{\tau}{N\rho_j}) - \lambda(y,0) - \frac{\tau}{N\rho_j} v(\lambda(x,\frac{\tau}{N\rho_j}),\frac{\tau}{N\rho_j})
&= a(\frac{i}{N}) - a(\frac{j}{N})+ O(\frac{1}{N^2}) \\
&=  a'(\frac{j}{N}) \frac{i-j}{N} = -  \frac{n}{N\rho_j} + O(\frac{1}{N^2}),
\end{align*}
 upon identifying $a'(\frac{j}{N}) = - 1/\rho_j$. It is now plain to check that the two formulas  \eqref{cauchy} and \eqref{sol-non-perturb} 
  match perfectly at the frontier between the mesoscopic and macroscopic regimes which corresponds to values of $t$ and $N$ such that
  $ 1/N \ll t \ll 1$.

\subsubsection{The stationary case}

We now consider an interesting special case where the bulk of the spectrum of $A$ 
has a Wigner semicircle density 
\begin{align*}
\rho_A(\lambda) = \frac{1}{2\pi} \sqrt{4-\lambda^2}\,. 
\end{align*}
We can slightly modify the definition of the matrix $X_t$ defined in \eqref{def-X}
 so that the limiting density of the eigenvalues is a Wigner semicircle at all time, $\rho(\lambda,t) = \frac{1}{2\pi} \sqrt{4-\lambda^2}$. 
 This can be done for instance by defining $X_t$ as the solution of the Ornstein-Uhlenbeck equation 
 \begin{align*}
dX_t :=-\frac{1}{2} X_t \, dt + dH_t\,, \quad \mbox{ with } \quad X_0= A\,,
\end{align*}
and where $H$ is a Hermitian Brownian motion as defined in \eqref{hermitian-bm}. 

In this special stationary case, the solution $u(x,y,t)$ of the evolution equation 
\eqref{overlap-eq-non-perturb} can be computed explicitly as was done in \cite{bourgade-univ}. 
For completeness we recall this computation here and  
propose an alternative derivation of the main identity \eqref{main-identity} in the appendix \ref{integral-chebyshev}.  

It is easier to work with the function $w$ introduced in \eqref{def-w}. 
Eq. \eqref{overlap-eq-non-perturb} may be rewritten in terms of $w(\lambda,\mu,t)$ as 
\begin{align}\label{eq-continue-3rd}
\partial_t \, w = - \mathcal{A} \, w \quad \mbox{ with } \quad w(\lambda,\mu,0)= 2\pi 
\frac{\delta(\lambda-\mu)}{\sqrt{4-\mu^2}}
\end{align}
 where $ \mathcal{A}$ is the linear operator on the space of smooth functions $f:[-2,2]\to \R$, defined as 
 \begin{align*}
\mathcal{A} \, f (\lambda):=  P.V. \int_{-2}^2 \frac{f(\lambda)-f(\nu)}{(\lambda-\nu)^2}  \frac{1}{2\pi} \sqrt{4-\nu^2}   \, d\nu\,. 
\end{align*}
The operator $\mathcal{A}$ already appeared in \cite{bourgade-univ} in a different context. 

It is easy to check (see appendix \ref{integral-chebyshev}) 
that the operator $\mathcal{A}$ has real eigenvalues $\{n/2, n \in \N \}$ and admits an {\it orthonormal} basis of eigenfunctions 
$(f_n)_{n\in \N}$ in the Hilbert space  
\begin{align*}
\mathcal{H}:=\left\{ f : [-2,2] \to \R: \int_{-2}^2 f(\lambda)^2 \sqrt{4-\lambda^2} \, d\lambda < +\infty\right\}
\end{align*}
endowed with the scalar product 
\begin{align*}
\langle f,g\rangle_{\mathcal{H}}:= \frac{1}{2\pi} \int_{-2}^2 f(\lambda) g(\lambda) \sqrt{4-\lambda^2} \, d\lambda\,. 
\end{align*}
The eigenfunctions $(f_n)_{n\in \N}$ are the Chebyshev polynomials of the second kind. More precisely, 
\begin{align*}
f_n(\lambda) = U_n(\frac{\lambda}{2})
\end{align*} 
where $U_n$ is the $n$-th Chebyshev polynomial of the second kind. 
Those  polynomials are the
orthogonal polynomials associated to the scalar $\langle \cdot, \cdot\rangle_{\mathcal{H}}$. 

\begin{figure}\center
\includegraphics[scale=0.8]{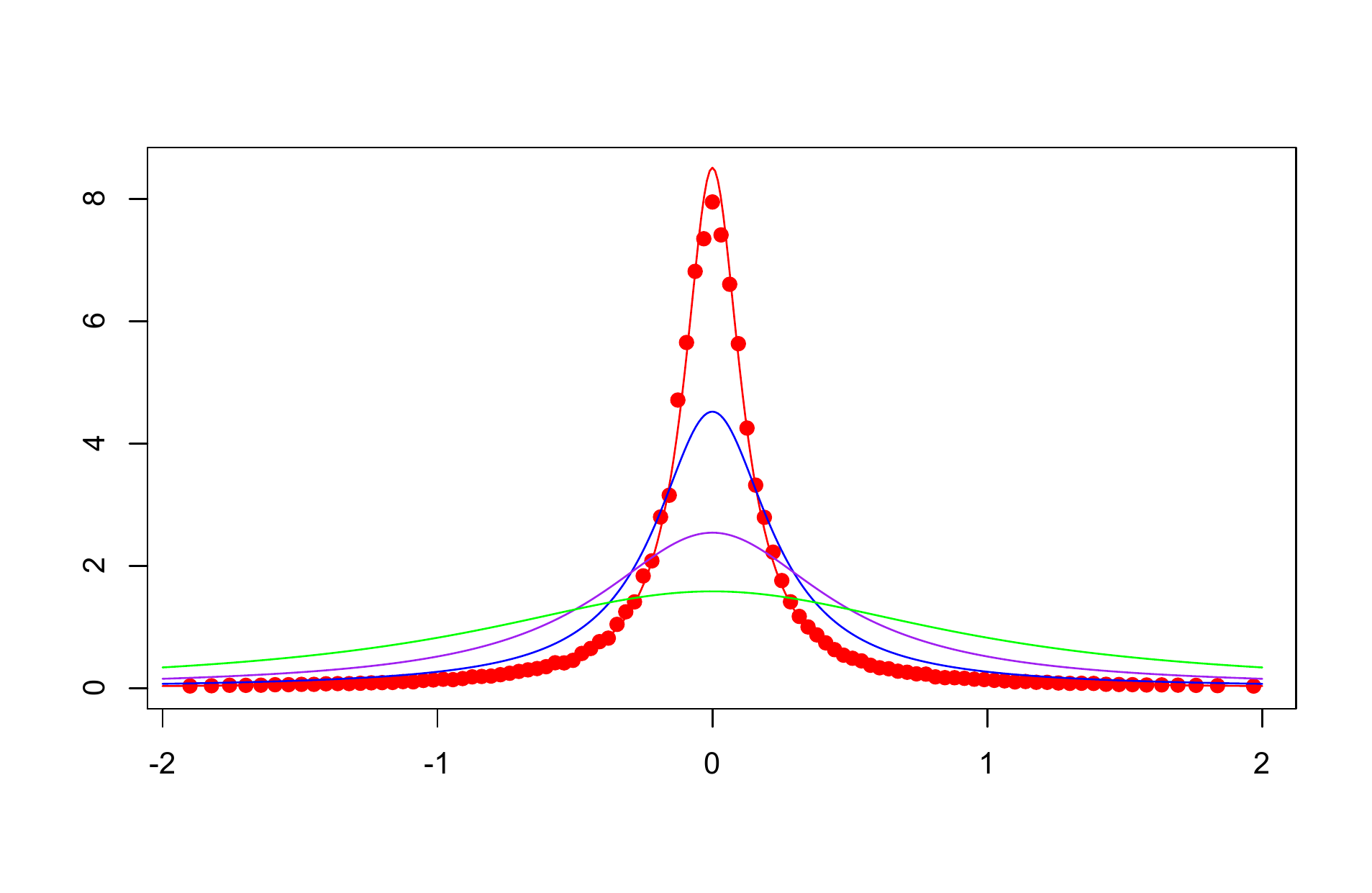}
\caption{(Color online). The plain curves represent the function $K_t(\lambda,0)$ as a function of $\lambda$ for $\mu=0$ 
for different values of $t=1/8$ (red), $t=1/4$ (blue), $t=1/2$ (purple), $t=1$ (green). 
One sees that the dispersion of the probability density $K_t(\lambda,0)$ 
increases as $t$ grows. The red points are computed using numerical simulations and represent
the rescaled overlaps $N \left[ \langle \psi_i^t |\psi_j^0\rangle^2 \right]$ as a function of
$\lambda_i\in[-2,2],$ $i=1,\dots,N$ and for $j=N/2$, $t=1/8$, $N=200$. 
The agreement with the corresponding theoretical curve as predicted in \eqref{conv-overlaps-continuous} is excellent. }\label{spreading-overlaps}
\end{figure}

The operator $\mathcal{A}$ is therefore self-adjoint in the Hilbert space $\mathcal{H}$
and we can easily compute the ``heat kernel'' $K_t(\lambda,\mu)$ of the continuous equation \eqref{eq-continue-3rd}
\begin{align*}
K_t(\lambda,\mu) = \sum_{n= 0}^{+\infty} \exp(-\frac{n}{2} t) U_n(\frac{\lambda}{2})U_n(\frac{\mu}{2})\,.
\end{align*} 
We can compute the sum of this latter series explicitly in terms of elementary functions 
\begin{align*}
K_t(\lambda,\mu) = \frac{1-e^{-t}}{1 - e^{-t/2} \lambda \mu + e^{-t}(\lambda^2+\mu^2-2) - 
\lambda \mu e^{-3t/2} + e^{-2t}} \,.
\end{align*}
The heat kernel $K_t(\lambda,\mu)$ also appeared in \cite[Page 462]{biane} in a different context. 
The solution of \eqref{eq-continue-3rd} is
\begin{align*}
w(\lambda,\mu,t) = \exp(- t \mathcal{A}) \frac{\delta(\lambda-\mu)}{\rho_A(\mu)} = K_t(\lambda,\mu) \,. 
\end{align*}
This finally gives a nice convergence result for the overlaps towards an explicit limit when $N\to \infty$ with $t$ fixed
independently of $N$,
\begin{align}\label{conv-overlaps-continuous}
N\,  \left[ \langle \psi_i^t |\psi_j^0\rangle^2 \right] \to K_t(\lambda,\mu) 
\end{align}
where  $\lambda:=\lim a(i_N/N)$ and $\mu:= \lim a(j_N/N)$. 
This formula was checked numerically for different values of $t$ and $N$ (see Fig. \ref{spreading-overlaps}).

\section{Isolated eigenvectors}\label{isolated}
In this section, we study the projections of a given initial eigenstate $|\psi_j^0\rangle$, associated to 
an eigenvalue $a_j$ lying outside the bulk of the spectrum of the initial matrix $A$, on the perturbed eigenvectors in 
the limit of large dimension $N$. 
If the eigenvalues $a_i$ are indexed in non-increasing order, then $j$ is a fixed finite integer which does not depend on $N$. 
To fix ideas and simplify notations, we will suppose that $j =1$: the eigenvalue $a_1$ is 
the largest spike (see Fig. \ref{fig-triangle}) of the matrix $A$.
Again we work with a sequence of matrices $A:=(A_N)_{N\in \N}$ such that Hypothesis \ref{hypothesis} holds.

\subsection{Principal component}\label{principal-component}
From the eigenvector evolution equation \eqref{evolution-psi_i}, we easily check that 
\begin{align*}
d \left[ \langle \psi_1^t|\psi_1^0\rangle\right] = -\frac{1}{2N} \sum_{k \neq 1} \frac{dt}{(\lambda_1(t)-\lambda_k(t))^2}
\left[ \langle \psi_1^t|\psi_1^0\rangle\right]\,.
\end{align*}
This ordinary differential equation is easily solved and, using the initial condition,  
we obtain the following equality, valid for any finite $N$,
\begin{align*}
\left[ \langle \psi_1^t|\psi_1^0\rangle\right] = \exp\left(-\frac{1}{2N} \int_0^t \sum_{k \neq 1} \frac{ds}{(\lambda_1(s)-\lambda_k(s))^2}\right)\,. 
\end{align*}
Sending $N\to \infty$, we easily get, using the results explained in subsection \ref{density-trajectory}, the almost sure (with respect to the eigenvalues) convergence 
\begin{align}\label{conv-vector-spike}
\left[ \langle \psi_1^t|\psi_1^0\rangle\right] \to 
\exp\left(-\frac{1}{2} \int_0^t ds \int_\R \frac{\rho(\lambda,s)}{(\lambda_1(s)-\lambda)^2} \, d\lambda \right)
\end{align}
where $(\lambda_1(s))_{0\le s \le t}$ is the limiting trajectory of the first spike (already described in \eqref{trajectory-spikes}) 
such that 
\begin{align}\label{spike1}
\dot{\lambda}_1(s) = \int_\R \frac{\rho(\lambda,s)}{\lambda_1(s)-\lambda} d\lambda\,, \quad \lambda_1(0)=a_1\,. 
\end{align}
The convergence \eqref{conv-vector-spike} can of course be extended from $j=1$ to any finite fixed value of $j$, with 
similar asymptotic formulas.  

%Even though $t$ is of order $1$ so that the perturbation $H_t$ has a macroscopic size, 

We see that, if $t$ is small enough so that the spike is still isolated from the bulk at time $t$, then 
 the overlap between the initial top eigenvector and its perturbed version
 does not vanish in the large $N$ limit even though $t$ and $H_t$ have macroscopic sizes, 
in contrast with the bulk overlaps which were of order $1/N$. 

\subsection{Transverse components} \label{transverse-components}
We now consider the overlaps between the initial top eigenvector $|\psi_1^0\rangle$ and the perturbed eigenvectors 
$|\psi_i^t\rangle$ for $i\neq 1$. The eigenvalue $\lambda_1(t)$ is isolated from the other eigenvalues
so that we expect the overlaps between the corresponding perturbed and non-perturbed eigenvectors to be 
microscopic of order $1/N$.  

To see this, we start again from the overlap equation \eqref{eq-overlaps}.
As before, we set $w(\lambda,t)=u(x,t):= \lim_{N\to \infty} \langle \langle \psi_1^0 | \psi_i^t\rangle^2 \rangle$ where 
$x:= \lim i_N/N$
 ($\lambda:=\lambda(x,t)$ is in the continuous part of the density $\rho(\cdot,t)$ and is the limit point of the eigenvalue $\lambda_i(t)$).  
 We also denote by $f(t)$ the limit of $\langle \langle \psi_1^0 | \psi_1^t\rangle^2 \rangle$ when $N\to +\infty$. 
 It is easy to compute $f(t)$ thanks to the overlap equation
 \begin{align}\label{def-f}
f(t) := \exp\left(- \int_0^t ds \int_\R \frac{\rho(\lambda,s)}{(\lambda_1(s)-\lambda)^2} \, d\lambda \right)\,.  
\end{align}
This convergence together with \eqref{conv-vector-spike} imply 
that $ \langle \psi_1^0 | \psi_1^t\rangle$ converges almost surely towards $\sqrt{f(t)}$ when $N\to +\infty$. 

We now easily derive the  Cauchy problem satisfied by the limiting 
family of overlaps  $u(x,t)$: for any $x\in [0,1]$,  
\begin{align}\label{cauchy-pb-transverse}
\partial_t \, u(x,t) = P. V. \int_{0}^1 \frac{u(y,t)-u(x,t)}{(\lambda(y,t)-\lambda(x,t))^2} \, dy + \frac{f(t)}{(\lambda_1(t)-\lambda(x,t))^2}\,, \quad
u(x,0) = 0\,. 
\end{align}
Note that the solution of \eqref{cauchy-pb-transverse} satisfies $u(x,t) \ge 0$ for all $t\ge 0$ and any $x$ 
in the bulk of the spectrum, as it should be for a mean squared overlap. 

We have the following almost sure convergence of the overlaps 
\begin{align}\label{conv-transverse}
N \left[ \langle \psi_1^0 | \psi_i^t \rangle^2 \right] \to u(x,t)
\end{align} 
where $x:=\lim i_N/N$  as $N\to+ \infty$. 
This result \eqref{conv-transverse} is equivalent to the almost sure weak convergence when $N\to \infty$ of the 
probability measure on the interval $[0,1]$ 
\begin{align}\label{prob-projections-transv}
N \sum_{i \neq 1} \left[ \langle \psi_i^t|\psi_1^0\rangle^2\right] \delta_{i/N}(dx) \Rightarrow u(x,t)\, dx
\end{align}
where $u(\cdot,t)$ is the (unique) solution to the Cauchy problem \eqref{cauchy-pb-transverse}.

\subsection{Gaussian fluctuations of the principal component}\label{gaussian-fluctuations}
Using the convergence \eqref{conv-transverse} of the transverse overlaps, we can compute the higher order moments of the 
principal component  and deduce that 
the random variable $\langle \psi_1^t|\psi_1^0\rangle$ is asymptotically a Gaussian variable with mean value $\sqrt{f(t)}$ as defined in 
\eqref{def-f}
and variance of order $1/N$ that we are able to compute explicitly.  

In this subsection , we work with a time $t>0$ small enough so  
that the spike $(\lambda_i(s))_{0\le s\le t}$ has not yet been swallowed by the limiting 
bulk density $(\rho(\lambda,s))_{0\le s\le t}$ of 
the Gaussian matrix process $(X_s)_{0 \le s\le t }$. This critical time $t_c$ was explicitly computed in section \ref{factor-model}
in the case of a small initial rank for the matrix $A$.
%In the following subsection, we compute explicitly the 
%time at which the largest eigenvalue is caught back by the bulk of the spectrum in the particular case 
%where the matrix $A$ has a small rank. 

For such a time $t<t_c$, we shall now prove that, {\it almost surely with respect to the eigenvalues path} $(\lambda_i(s))_{s<t},i=1,\dots,N$, the random variable 
\begin{align*}
\sqrt{N} \, \left( \langle \psi_1^t | \psi_1^0 \rangle -\left[ \langle \psi_1^t | \psi_1^0 \rangle \right] \right)
\end{align*}
converges weakly towards a centered Gaussian distribution with variance 
\begin{align*}
g^{(2)}(t) := \int_0^t ds \exp\left(- \int_s^t \int_\R \frac{\rho(\lambda,u)}{(\lambda_j(u)-\lambda)^2} \, d\lambda\, du \right)  \int_\R \frac{w(\mu,s)}{(\lambda_j(s)-\mu)^2} \rho(\mu,s) d\mu\,,
\end{align*}
where \begin{itemize}
\item $(\lambda_1(s))_{0\le s \le t}$ is the limiting trajectory of the largest eigenvalue satisfying \eqref{trajectory-spikes}; 
\item $(\rho(\lambda,s))_{0\le s \le t,\lambda\in \R}$ is the limiting bulk density trajectory satisfying \eqref{burgers};
\item $(u(x,s))_{0\le s \le t,x\in [0,1]}=(w(\lambda,t)_{0\le s\le t, \lambda\in \R})$ 
is the function describing the limiting transverse overlaps satisfying the evolution equation \eqref{cauchy-pb-transverse}.
\end{itemize}

 If we denote by 
\begin{align*}
g_N^{(2)}(t) :=  \left[ \left (\langle \psi_1^t | \psi_1^0 \rangle -\left[ \langle \psi_1^t | \psi_1^0 \rangle \right] \right)^2 \right], 
\end{align*}
we can easily check with the It\^o's formula that 
\begin{align}\label{eq-diff-g}
d g_N^{(2)}(t) = - \frac{1}{N} \sum_{k \neq 1} \frac{dt}{(\lambda_1-\lambda_k)^2} g_N^{(2)}(t) + \frac{1}{N}  h_N(t) \, dt
\end{align}
where 
\begin{align*}
h_N(t):= \sum_{k \neq 1} \frac{[ \langle \psi_k^t |\psi_1^0 \rangle^2 ]}{(\lambda_1-\lambda_k)^2}\,. 
\end{align*}
It is straightforward to solve the ordinary differential equation \eqref{eq-diff-g} 
\begin{align}
g_N^{(2)}(t) = \frac{1}{N} \int_0^t \exp\left(-\frac{1}{N} \int_s^t \sum_{k \neq 1} \frac{du}{(\lambda_1(u)-\lambda_k(u))^2} \right) h_N(s) \, ds \,.
\label{sol-g_N}
\end{align}
The limit of the function $h_N(s)$ is easily computed thanks to the results obtained in the previous subsection. We have the following almost sure 
convergence (with respect to the eigenvalues)
\begin{align*}
h_N(s) \to h(s):=  \int_\R \frac{w(\lambda,s)}{(\lambda_1(s)-\lambda)^2} \rho(\lambda,s) d\lambda
\end{align*}
where $u$ satisfies the Cauchy problem \eqref{cauchy-pb-transverse} and $\lambda_1$ is the limiting trajectory of the first spike. 
From \eqref{sol-g_N}, we deduce the almost sure convergence of the rescaled function $N \, g_N^{(2)}(t)$ when $N\to \infty$,
\begin{align}\label{def-g2}
N \, g_N^{(2)}(t) \longrightarrow g^{(2)}(t):= \int_0^t \exp\left(- \int_s^t \int_\R \frac{\rho(\lambda,u)}{(\lambda_1(u)-\lambda)^2} \, d\lambda\, du \right) h(s) \, ds \,. 
\end{align}
With a similar method, we can check that the $n$-th moment 
\begin{align*}
g_N^{(n)}(t) := \left[ \left (\langle \psi_1^t | \psi_1^0 \rangle -\left[ \langle \psi_1^t | \psi_1^0 \rangle \right] \right)^n \right]
\end{align*}
satisfies 
\begin{align}\label{eq-diff-g-n}
d g_N^{(n)}(t) = - \frac{n}{2N} \sum_{k \neq 1} \frac{dt}{(\lambda_1-\lambda_k)^2} g_N^{(n)}(t) + \frac{n(n-1)}{2N}\, g_N^{(n-2)}(t) \, 
 h_N(t) \, dt\,. 
\end{align}
The ordinary differential equation \eqref{eq-diff-g-n} can be solved and we get for any $t\ge 0$ and $n \in \N$, 
\begin{align}\label{recursion-relation}
g_N^{(n)}(t) = \frac{n(n-1)}{2N} 
\int_0^t \exp\left(-\frac{n}{2N} \int_s^t \sum_{k \neq 1} \frac{du}{(\lambda_1(u)-\lambda_k(u))^2} \right) g_N^{(n-2)}(s) \, 
 h_N(s) \, ds \,. 
\end{align}
Using this recursion relation \eqref{recursion-relation}, we can now prove iteratively on $n\in \N$ that, as $N\to \infty$,
\begin{align*}
N^{n/2} g_N^{(n)}(t) \longrightarrow g^{(n)}(t)
\end{align*}
%almost surely 
%with respect to the eigenvalues trajectories and
where the sequence $(g^{(n)}(t))_{n \in \N}$ satisfies the recursion relation 
 \begin{align}\label{limiting-recursion-relation}
g^{(n)}(t) = \frac{n(n-1)}{2} 
\int_0^t \exp\left(- \frac{n}{2}\int_s^t \int_\R \frac{\rho(\lambda,u)}{(\lambda_1(u)-\lambda)^2} d\lambda du \right) g^{(n-2)}(s) \, 
 h(s) \, ds \,. 
\end{align}
To prove that the limiting distribution of the random variable 
\begin{align*}
\sqrt{N} \left( \langle \psi_1^t | \psi_1^0 \rangle -\left[ \langle \psi_1^t | \psi_1^0 \rangle \right]  \right)
\end{align*} 
is indeed Gaussian with mean $0$ and variance $g^{(2)}(t)$, 
it suffices to check that the sequence $(g^{(n)}(t))_{n \in \N}$ corresponds to 
the moments of a Gaussian variable with mean $0$ and variance $g^{(2)}(t)$, i.e. that 
for all $n\in \N$, 
\begin{align}\label{moments-gaussian}
&g^{(2n+1)}(t) = 0 \notag \,, \\
&g^{(2n)}(t) = (2n-1) \, g^{(2)}(t) \, g^{(2n-2)}(t) \,. 
\end{align}
It is straightforward to see that the distribution is symmetric with zero odd moments iteratively using the recursion relation 
\eqref{recursion-relation} valid for finite values of $N$ and the
convergence to $0$ of the first moment when $N\to \infty$ obtained in Eq. \eqref{conv-vector-spike} subsection \ref{principal-component}. 

The proof that the even moments satisfy the relation \eqref{moments-gaussian} for any $n$ can be found in Appendix \ref{proof-even-moments}. 
 
An alternative proof can be done {\it ad hoc} with the characteristic function 
\begin{align*}
F_N(\xi,t) = \left[\exp\left({\rm i} \xi \sqrt{N} \left( \langle \psi_1^t |\psi_1^0 - \left[\langle \psi_1^t |\psi_1^0\rangle   \right]  \right) \right)\right].
\end{align*} 
 It is plain to check thanks to It\^o's formula that the function $F_N$ satisfies the partial differential equation
 \begin{align*}
\frac{\partial }{\partial t} F_N(\xi,t) = -\frac{\xi}{2N}  \frac{\partial }{\partial \xi} F_N(\xi,t) \sum_{k\neq 1} \frac{1}{(\lambda_1-\lambda_k)^2} 
- \frac{\xi^2}{2} h_N(t) F_N(\xi,t)\,. 
\end{align*}
In the scaling limit $N\to \infty$, this equation becomes 
\begin{align*}
\frac{\partial }{\partial t} F(\xi,t) = -\frac{\xi}{2}  \frac{\partial }{\partial \xi} F(\xi,t)\int_\R  \frac{\rho(\lambda,t)}{(\lambda_1(t)-\lambda)^2} d\lambda 
- \frac{\xi^2}{2} h(t) F(\xi,t)
\end{align*} 
 which is clearly satisfied by the Gaussian characteristic function $F(\xi,t) = \exp(-\frac{\xi^2}{2}g^{(2)}(t))$. 

\subsection{Estimation of the main factors}\label{factor-model-vectors}
As an illustration of the results obtained in the 
previous subsection, we come back on the factor model. 
In section \ref{factor-model}, we have seen that the limiting density of eigenvalues is the Wigner semicircle 
with radius $2\sqrt{t}$ at time $t$ and that the limiting trajectories of the spikes are $\lambda_j(t)=a_j+t/a_j$. 

It turns out that the limiting mean square overlap between the first non-perturbed and perturbed eigenvectors (respectively 
$|\psi_1^0\rangle$ and $|\psi_1^t\rangle$) can also be computed analytically. 

From \eqref{conv-vector-spike}, we obtain 
\begin{align*}
\left[ \langle \psi_1^t|\psi_1^0\rangle\right] \to 
\exp\left(-\frac{1}{2} \int_0^t \frac{ds}{2 \pi s} \int_\R \frac{\sqrt{4s-\lambda^2}}{(a +\frac{s}{a}-\lambda)^2} \, d\lambda \right)
= \sqrt{\max(1-\frac{t}{a_1^2},0)} \,.  
\end{align*}
We see that the information contained in the perturbed eigenvector is completely lost at the time $t_c=a_1^2$ 
when the spike $\lambda_1(t_c)$ is
swallowed by the Wigner sea. 

From the results of subsection \ref{gaussian-fluctuations}, we know that the random variable 
$\sqrt{N} \langle \psi_1^0 |\psi_1^t\rangle$ has Gaussian fluctuations in the large $N$-limit around its mean asymptotic value $1-t/a_1^{2}$ 
for $t\le a_1$.

 \begin{figure}\center
\includegraphics[scale=0.95]{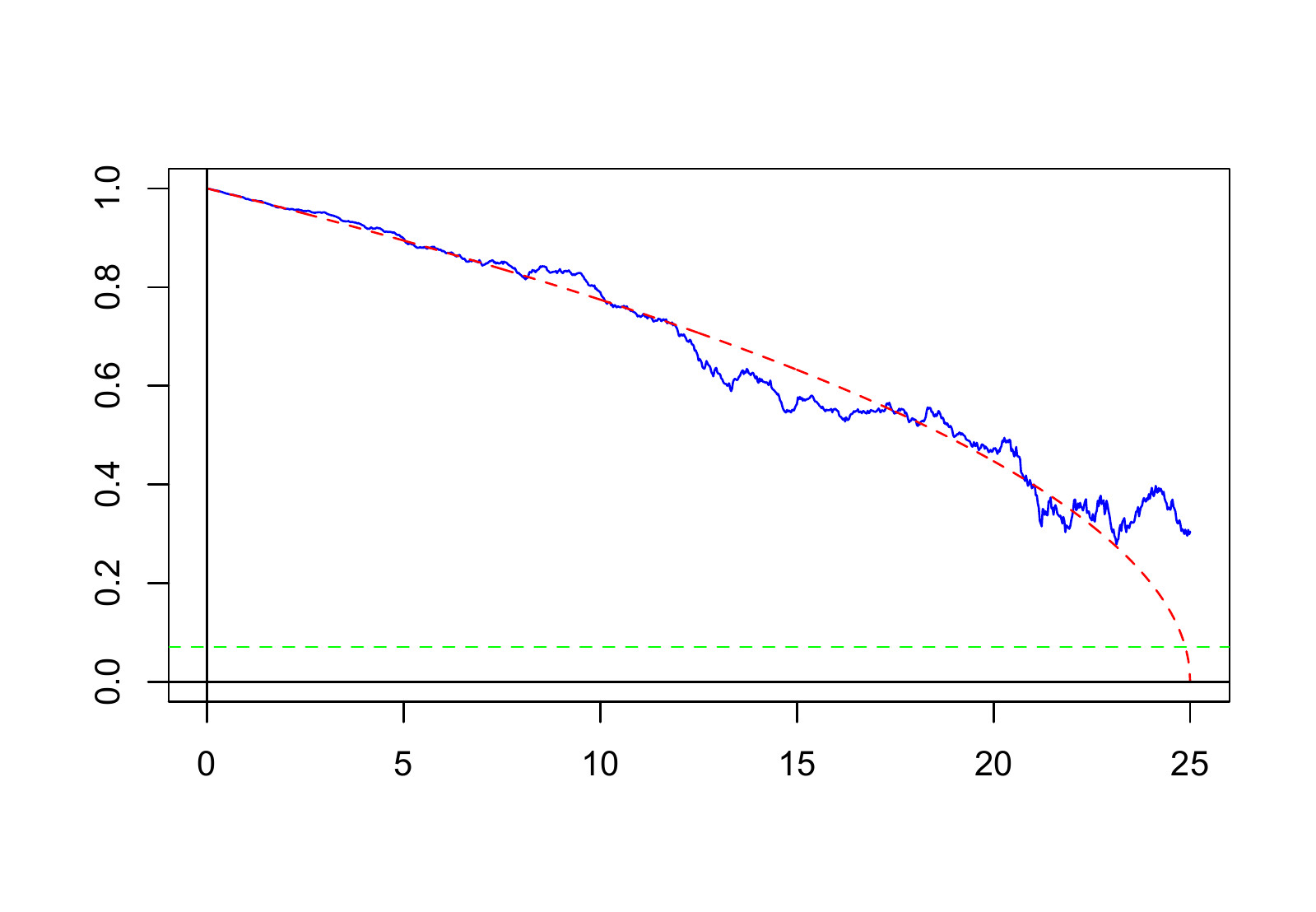}
\caption{(Color online). 
Numerical simulation of the process $\langle\psi_1^t|\psi_1^0\rangle$ as a function of time $t\in [0,a_1^2]$ (blue line)
together with the theoretical limiting curve $\sqrt{1-t/a_1^2}$ for $t \leq t_c=25$ (red dashed line). 
The matrix $A$ has only one non zero eigenvalue $a_1=5$ and 
the dimension is $N=200$. The horizontal green dashed line is $1/\sqrt{N}$.
The agreement is good away from the right end point $t_c=a_1^2=25$ where a phase transition must occur near the critical point. 
We see that the convergence holds almost surely as predicted in subsection \eqref{transverse-components} 
from the convergence of the second moment. }\label{isolated-ev-simul}
\end{figure}

The probability measure on $[0,1]$ with weights given by the overlaps of $|\psi_1^0\rangle$ (introduced in \eqref{prob-projections-transv}) 
converges weakly almost surely 
\begin{align*}
N \sum_{i \neq 1} \left[ \langle \psi_1^0 |\psi_i^t\rangle^2 \right] \delta_{i/N}(dx) \Rightarrow u(x,t)\, dx 
\end{align*}
where the function $u(x,t)$ (which contains the information 
the microscopic transverse components of the initial vector $|\psi_1^0\rangle$) satisfies, 
for any $t\ge 0$ and $x\in [-2\sqrt{t},2\sqrt{t}]$, 
\begin{align*}
\partial_t \, u(x,t) &=  P. V. \int\limits_{0}^{1} \frac{u(y,t)-u(x,t)}{(\lambda(y,t)-\lambda(x,t))^2} \, dy 
+ \frac{\max(1-\frac{t}{a_1^2},0)}{(a_1+ \frac{t}{a_1} - \lambda(x,t))^2} 
 \,, \\ 
u(\cdot,0) &= 0  \,,  
\end{align*}
where the quantile function $\lambda(x,t)$ satisfies for any $x\in (0,1)$, 
\begin{align*}
x = \frac{1}{2\pi t} \int\limits_{\lambda(x,t)}^{2\sqrt{t}} \sqrt{4t-\lambda^2} \, d\lambda\,.
\end{align*}
%To solve this evolution equation, we first consider the homogeneous equation whose general solution $u_0(\cdot,t)$
%can be computed explicitly: for  
%$t\ge 0$ and $x\in [-2\sqrt{t},2\sqrt{t}]$, 
%\begin{align*}
%u_0(x,t)= \sum_{n = 0}^{+\infty} \exp(-\frac{n}{2} t)  U_n(0) U_n(\frac{x}{2\sqrt{t}})
%\end{align*}
%Setting $v(x,t):= u(x\sqrt{t},t)$, $\partial_x v(x,t) =\sqrt{t}\,  \partial_x u(x\sqrt{t},t),$ 
%$\partial_t v(x,t) = \frac{x}{2\sqrt{t}} \partial_x u(x\sqrt{t},t) + \partial_t u(x\sqrt{t},t)$  
%Eq. \eqref{cauchy-pb-transverse-2} may be rewritten in terms of $v$ as  
%\begin{align}\label{cauchy-pb-transverse-v}
%\partial_t &\, v(x,t) - \frac{x}{2t} \partial_x v(x,t) \notag
%\\ &=  \frac{1}{2\pi } P. V. \int_{-2}^{2} \frac{v(y,t)-v(x,t)}{(y-x)^2} \sqrt{4-y^2} dy 
%+ \frac{1}{(a_1+ \frac{t}{a_1} - x\sqrt{t})^2} 
%\sqrt{\max(1-\frac{t}{a_1^2},0)}  \\ 
%v(\cdot,0) &= 0 \notag \,. 
%\end{align}

\section{Extension to Covariance matrices and conclusion}\label{covariance}
We now consider a similar problem motivated by applications in finance. We want to estimate 
a $N\times N$ positively definite matrix $C$ (covariance matrix) 
from the observation of a sequence of independent and identically distributed centered Gaussian vectors 
$r^t:=(r_1^t,r_2^t,\dots,r_N^t),t\in \N$ whose covariance matrix is $C$, i.e. such that for any $1\le i,j\le N$ and $t\in \N$, 
\begin{align*}
\E[r_i^t r_j^t] = C(ij)\,. 
\end{align*}
In the context of finance, $N$ is the number of stocks in the financial market under consideration and 
 the Gaussian variable $r_i^t$ is the return of the $i$-th stock on the $t$-th day. The covariance matrix $C$ is unknown 
 and an important issue for risk control is to estimate it from the empirical data. Without loss of generality, we can suppose 
 that the matrix $C$ is diagonal and denote by $ c_1\ge c_2\ge\dots \ge c_N\ge 0$ its real non-negative eigenvalues. 

The basic idea is to 
form the empirical covariance matrix $E_T$ from a sequence of $T$ observations of the returns, defined as 
\begin{align*}
E_T(ij) := \frac{1}{T} \sum_{t=1}^T r_i^t r_j^t\,. 
\end{align*}
If the number of stocks $N$ is fixed and if we have access 
to a very long sequence of observations of the returns of length  $T\to +\infty$, then $E_T$ gives a perfect estimation of   
$C$. Indeed, the central limit theorem implies that 
\begin{align*}
E_T(ij) \to_{T\to+ \infty} C(ij)\,. 
\end{align*}
 Nevertheless, in practical applications, one only has access to finite length datasets and the number of stocks (or variables) $N$
 is often quite large, comparable to the length $T$ of the time series. The relevant framework 
 to be considered is the case when both $N$ and $T$ tend to infinity but with a fixed ratio $q$, i.e.
 \begin{align}\label{double-limit}
N\to +\infty\,, \quad T\to +\infty \quad  \mbox{ with } \quad \frac{N}{T}\to q\,. 
\end{align} 
Let us exhibit the similarities of this present problem with the simpler one investigated in the previous sections. 
The matrix $E_T$ can be decomposed as follows 
\begin{align}\label{decomposition-E}
E_T(ij) = C(ij) + \mathcal{E}_T(ij) \quad \mbox{ with } \quad \mathcal{E}_T(ij):= \frac{1}{T} \sum_{t=1}^T r_i^t r_j^t - C(ij) \,. 
\end{align}
 We easily check using the central limit theorem that 
 \begin{align*}
\mathcal{E}_T(ij) \sim \frac{1}{\sqrt{T}} = \sqrt{\frac{q}{N} }\,. 
\end{align*}
Therefore the matrix $E_T$ may be seen as an additive random perturbation of the deterministic matrix $C$ 
as in the previous case \eqref{def-X}. The time parameter $t$ is now replaced by the {\it quality factor} $q$. 
The covariance structure  of the entries of $\mathcal{E}$ is nevertheless slightly more complicated in the present case (the entries of the perturbation $H_t$ in \eqref{def-X} are independent and identically distributed). 
A straightforward computation leads to
\begin{align*}
\E[\mathcal{E}(ij) \mathcal{E}(kl)] = q \frac{C(ik) C(jl) + C(il) C(jk)}{N}\,. 
\end{align*}
We can deduce from the latter formula that the variables $\mathcal{E}_T(ij) $ and $\mathcal{E}_T(kl)$ 
are (asymptotically in the double limit \eqref{double-limit}) independent whenever the couple $(i,j)$ is different from $(k,l)$ and from $(l,k)$.   
If $(i,j)=(k,l)$ or $(i,j)=(l,k)$ with $i\neq j$, the variance is 
\begin{align*}
\E[\mathcal{E}_T(ij)^2] = q \frac{c_i c_j}{N}\,. 
\end{align*}
If $i=j$, 
\begin{align*}
\E[\mathcal{E}_T(ii)^2] = 2 q \frac{c_i^2}{N}\,. 
\end{align*}
We conclude that the structure of covariance of the entries is very similar to the previous case $X_t=A+H_t$ 
(the entries of $\mathcal{E}$ are independent) but with a different non-homogeneous variance profile 
for the entries of $\mathcal{E}$. 

The question is now whether one can construct a  Hermitian matrix process  going from $C$ to $E_{ij}(T)$. 
As mentioned above, the time is $q$ and we will now write $E_q$ instead of $E_T$. 
At time $q=0$ corresponding to $T=+\infty$, 
\begin{align}\label{final-cond}
E_0= C\,. 
\end{align}
The idea is now to see the matrix $E_{q-dq}$ as an additive perturbation of the matrix $E_q$.
Note that we do not proceed through the forward way but rather backward and the condition \eqref{final-cond} at $q=0$  
should be seen as a final condition. The initial condition in $q=+\infty$ is
\begin{align*}
E_{+\infty} = 0
\end{align*}
and the dynamical evolution goes backward in time $q$. 
 From \eqref{decomposition-E}, it is easy to see that, in the double scaling limit \eqref{double-limit}, we have the 
 following evolution  
 \begin{align}\label{evolution-E}
E_{q-dq} -E_{q} =  (C-E_q ) \,   dq  + \frac{1}{\sqrt{N}} dG_q\,,
\end{align}
where $(G_q)_{q\ge 0}$ is a Hermitian matrix process whose entries are independent Brownian motions 
with quadratic variations 
\begin{align*}
\E[ dG_q(ij)^2 ] =  \frac{c_i c_j}{N} dq \quad \mbox{ if } \quad i \neq j\,, \\ 
\E[ dG_q(ii)^2 ] = 2 \frac{c_i^2}{N} dq \quad \mbox{ otherwise}\,.
\end{align*}
 The empirical matrix $E_q$ performs a Ornstein-Uhlenbeck process around its asymptotic value $C$ reached in the limit 
 $q\to 0$. At the matrix level, the process is explicit and one could expect to take advantages of this 
 simple description for the empirical matrix process  $(E_q)_{q\in \R_+}$. Nevertheless, 
 the evolution \eqref{evolution-E} is not isotropic if $C$ is different from the identity matrix $I$
 as the empirical matrix is pushed in the direction of the matrix $C$. 
 If $C\neq I$, this non-isotropy implies that the eigenvalues evolution equation is not autonomous and 
 depends also on the eigenvector process.

Nevertheless, this new description of the empirical process $E_q$ could be used to recover recent results obtained on 
the relationship between the sample eigenvalues/eigenvectors of the empirical matrix
$E_q$ with the population eigenvalues/eigenvectors of the matrix $C$ in \cite{sandrine,paul,benaych,bbp,bloemendal-1,bloemendal-2}. 

A relevant choice for $C$ used in practical applications 
is a factor model where one supposes
that the spectrum of  $C$ contains a finite number of spikes $c_1\ge \dots \ge c_k \ge 1$
together with an eigenvalue $1-\delta$ \footnote{The parameter $\delta$ is usually chosen such 
that the trace of the covariance matrix 
$C$ is $N$ i.e. such that $(1-\delta)(N-k) + \sum_{i=1}^k c_i = N$.} with multiplicity $N-k$. The  spikes $c_1 \ge c_2  \dots \ge c_k$ 
are respectively associated
to eigenvectors $|\psi_1^0 \rangle, |\psi_2^0\rangle, \dots,  |\psi_k^0\rangle$ where $|\psi_1^0\rangle$ corresponds to 
the market mode while the $ |\psi_i^0\rangle, i=2, \dots,k$ usually contain the information on 
the economic sectors in the context of finance \cite{review-jpb}.
The statistics of the empirical estimation of the spikes were first investigated in \cite{bbp} whose study was later completed in 
\cite{bloemendal-1,bloemendal-2}. The overlaps between the sample and population eigenvectors 
associated to eigenvalues lying in the bulk of the continuous part of the empirical  
spectrum (whose limiting shape is the Marchenko-Pastur density) were computed in \cite{sandrine}. 
The overlap between the population and sample top eigenvectors associated to the spikes were investigated in 
\cite{paul,benaych}. 

We think it would be very interesting to recover those results following the lines of our present approach 
based on the evolution equation 
\eqref{evolution-E}.  
We leave this challenging problem for future research.

%Extension to covariance matrices. 
% Model with factors: principal component in \cite{paul,benaych}, no transverse components can be defined in this case.
% General covariance spiked matrix $C$ with continuous part, MP formula for general $C$, same picture as described above
%  but no convenient brownian motions in this case. Mention that $\tau=N/T$ may be seen as time, and that one has 
%  a Ornstein Uhlenbeck process of $E^\tau$ around $C$.  

\appendix

\section{Proof of the local resolvent evolution \eqref{burgers-local}} \label{burgers-local-proof}
Thanks to It\^o's formula, we get
\begin{align*}
\partial _t U_N(z,t) &= \frac{1}{N} \sum_{i=1}^N \frac{1}{z-\lambda_i(t)} \sum_{k \neq i} \frac{u_{k|j}(t)-u_{i|j}(t)}{(\lambda_k-\lambda_i)^2}  + 
\sum_{i=1}^N \frac{u_{i|j}(t)}{(z-\lambda_i(t))^2}\frac{d \lambda_i}{dt} + \frac{2}{\beta N} \sum_{i=1}^N \frac{u_{i|j}(t)}{(z-\lambda_i(t))^3}  \\
&= \frac{1}{2N} \sum_{i\neq k} \frac{u_{k|j}(t)-u_{i|j}(t)}{\lambda_i-\lambda_k} \frac{1}{(z-\lambda_i)(z-\lambda_k)}\\
&+ \sum_{i=1}^N \frac{u_{i|j}(t)}{(z-\lambda_i(t))^2} \left( \sqrt{\frac{2}{\beta N}} \frac{dB_i}{dt} + \frac{1}{N} \sum_{k \neq i} \frac{1}{\lambda_i- \lambda_k}\right) + \frac{1}{\beta N} \partial_z^2 U_N(z,t)
\end{align*}
where we have used the classical symmetrization trick to obtain the second line. 

Now, the new trick is to rewrite the first term as 
\begin{align*}
\frac{1}{2N}& \sum_{i\neq k} \frac{u_{k|j}(t)-u_{i|j}(t)}{\lambda_i-\lambda_k} \frac{1}{(z-\lambda_i)(z-\lambda_k)}
\\&= \frac{1}{2N}  \sum_{k= 1}^N \frac{u_{k|j}(t)}{z-\lambda_k} \sum_{i \neq k} \frac{1}{(\lambda_i-\lambda_k)(z-\lambda_i)} - \frac{1}{2N}  \sum_{i= 1}^N \frac{u_{i|j}(t)}{z-\lambda_i} \sum_{k \neq i} \frac{1}{(\lambda_i-\lambda_k)(z-\lambda_k)} \,. 
\end{align*}
We notice that 
\begin{align*}
\sum_{i \neq k} \frac{1}{(\lambda_i-\lambda_k)(z-\lambda_i)}  = \frac{1}{z-\lambda_k} 
\sum_{i \neq k} \frac{1}{z-\lambda_i} + \frac{1}{\lambda_i-\lambda_k}
\end{align*}
and
\begin{align*}
\sum_{k \neq i} \frac{1}{(\lambda_i-\lambda_k)(z-\lambda_k)}  = -\frac{1}{z-\lambda_i} 
\sum_{k \neq i} \frac{1}{z-\lambda_k} - \frac{1}{\lambda_i-\lambda_k}\,. 
\end{align*}
Therefore we deduce that 
\begin{align*}
\frac{1}{2N} \sum_{i\neq k} &\frac{u_{k|j}(t)-u_{i|j}(t)}{\lambda_i-\lambda_k} \frac{1}{(z-\lambda_i)(z-\lambda_k)}
= \sum_{k=1}^N \frac{u_{k|j}(t)}{(z-\lambda_k)^2} \frac{1}{N} \sum_{i \neq k} \frac{1}{z-\lambda_i}
+  \sum_{k=1}^N \frac{u_{k|j}(t)}{(z-\lambda_k)^2}  \frac{1}{N} \sum_{i \neq k} \frac{1}{\lambda_i-\lambda_k}
\\
&= \sum_{k=1}^N \frac{u_{k|j}(t)}{(z-\lambda_k)^2}  \left(G_N(z,t) - \frac{1}{N} \frac{1}{z-\lambda_k}\right) 
+  \sum_{k=1}^N \frac{u_{k|j}(t)}{(z-\lambda_k)^2}  \frac{1}{N} \sum_{i \neq k} \frac{1}{\lambda_i-\lambda_k} 
\end{align*}
where 
$G_N(z,t) :=  \frac{1}{N} \sum_{i=1} \frac{1}{z-\lambda_i(t)}$.
As a conclusion, we deduce that 
\begin{align*}
\partial _t U_N(z,t) = - G_N(z,t) \partial_z U_N(z,t) + \sqrt{\frac{2}{\beta N}} \sum_{i=1}^N 
\frac{u_{i|j}(t)}{(z-\lambda_i)^2} \frac{dB_i}{dt} + \frac{1}{2N}(\frac{2}{\beta}- 1)\partial_z^2 U_N(z,t)\,. 
\end{align*}

\section{Integral relations for Chebyshev polynomials}\label{integral-chebyshev}
Denote by $U_n$ the $n$-th Chebyshev polynomial of the second kind such that 
$U_n(\cos \theta)=\frac{\sin (n+1) \theta}{\sin\theta}$ for $\theta\in \R$. 
For $z$ such that $\Im z >0$, we shall prove that, for any $n\in \N$, 
\begin{align*}
\mathcal{I}:= \frac{1}{2\pi} \int_{-2}^2 \frac{U_n(\frac{\mu}{2})}{z-\mu} \sqrt{4-\mu^2} d\mu  = \frac{n}{2} U_n(\frac{z}{2}) \,. 
\end{align*}
We have 
\begin{align*}
\mathcal{I}&= \frac{1}{\pi} \int_{0}^{2\pi} \frac{\sin \theta}{z-2\cos \theta} \sin((n+1)\theta)\, d\theta\\ 
&= \frac{1}{4\pi i} \oint_{|\zeta|=1} \frac{\zeta^2-1}{\zeta^{n+2} (\zeta^2 - z \zeta+1)} (\zeta^{2(n+1)}-1) d\zeta \,. 
\end{align*}
We can compute this integral thanks to the residues Theorem. We have three poles 
$p_0=0,p_1=(z+ \sqrt{z^2-4})/2, p_2=(z- \sqrt{z^2-4})/2$ where we choose the branch of the square root with positive imaginary part.  
Noting that $p_1p_2=1$, we can check that in fact $|p_1|>|p_2|$ and therefore $|p_2|<1$. 
We set $f(\zeta)=\frac{1}{1-z\zeta+\zeta^2}$. 
The generating function of the $U_n$ is
\begin{align*}
\sum_{k=0}^{+\infty} U_k(\frac{z}{2})\, \zeta^{k} = \frac{1}{1-z\zeta+\zeta^2}\,. 
\end{align*}
%It is easy to check that for any $k \in \N$, we have 
%\begin{align*}
%\frac{f^{(2k)}(0)}{(2k)!} = \sum_{j=0}^k \binom{k+j}{k-j} (-1)^{k-j} z^{2j}\,,\\
%\frac{f^{(2k+1)}(0)}{(2k+1)!} = \sum_{j=0}^k \binom{k+j+1}{k-j} (-1)^{k-j} z^{2j+1}\,. 
%\end{align*}
We now have to compute the three residues. The residue at $p_0$ is 
\begin{align*}
\mbox{Res}(p_0)= U_{n+1}(\frac{z}{2}) -U_{n-1}(\frac{z}{2}) = 2 \, T_{n+1}(\frac{z}{2}) 
\end{align*}
where $T_k$ is the $k$-th Chebyshev polynomial  of the first kind. 
The residue at $p_2$ is  
\begin{align*}
\mbox{Res}(p_2)&= \frac{1}{2^{n+1}} \left((z-\sqrt{z^2-4})^{n+1}- (z-\sqrt{z^2-4})^{-(n+1)} \right)\\
&=  \frac{1}{2^{n+1}} \left((z-\sqrt{z^2-4})^{n+1}- (z+\sqrt{z^2-4})^{n+1} \right)\\
&= - U_n(\frac{z}{2})\,  \sqrt{z^2-4} \\
&= - v_n(z)\, \sqrt{z^2-4}
 \,. 
\end{align*}
We finally obtain 
\begin{align*}
 \frac{1}{2\pi} \int_{-2}^2 \frac{U_n(\frac{\mu}{2})}{z-\mu} \sqrt{4-\mu^2} dy =   T_{n+1}(\frac{z}{2}) - \frac{1}{2}U_n(\frac{z}{2}) \, \sqrt{z^2-4} \,. 
\end{align*}
Differentiating with respect to $z$ and using the relation $T_{n+1}'=(n+1) U_n$, we get 
\begin{align*}
 \frac{1}{2\pi} \int_{-2}^2 \frac{U_n(\frac{\mu}{2})}{(z-\mu)^2} \sqrt{4-\mu^2} \, dy = -  \frac{n+1}{2} U_n(\frac{z}{2}) + \frac{1}{2} \left(\frac{1}{2}U_n'(\frac{z}{2})\, \sqrt{z^2-4} + U_n(\frac{z}{2}) \frac{z}{\sqrt{z^2-4}} \right) \,. 
\end{align*}
Similarly we have  
\begin{align*}
\frac{1}{2\pi} \int_{-2}^2 \frac{\sqrt{4-\mu^2}}{(z-\mu)^2} \,  d\mu = \frac{1}{2} (-1+ \frac{z}{\sqrt{4-z^2}})\,. 
\end{align*}
Therefore, 
\begin{align*}
\frac{1}{2\pi} \int_{-2}^2\frac{U_n(\frac{z}{2})-U_n(\frac{\mu}{2})}{(z-\mu)^2} \sqrt{4-\mu^2}\,  d\mu =\frac{n}{2} U_n(\frac{z}{2})  - \frac{1}{4} U_n'(\frac{z}{2}) 
\sqrt{z^2-4}\,. 
\end{align*}
Setting now $z=x+i\varepsilon$ and sending $\varepsilon \to 0$, we obtain keeping only the real part 
\begin{align}\label{main-identity}
\frac{1}{2\pi}P.V. \int_{-2}^2\frac{U_n(\frac{\lambda}{2})-U_n(\frac{\mu}{2})}{(\lambda-\mu)^2} \sqrt{4-\mu^2} \, d\mu 
=\frac{n}{2} U_n(\frac{\lambda}{2}) \,.   
\end{align}

\section{Proof of the Gaussian fluctuations announced in section \ref{gaussian-fluctuations}}\label{proof-even-moments}

For the even moments, we start by proving the relation \eqref{moments-gaussian} for $n=2$ i.e. that $g^{(4)}(t) = 3 g^{(2)}(t)^2$. 
To simplify notations in the following computations, we introduce the function $\varphi$ such that for $u \ge 0$, 
$\varphi(u) := \int_\R \frac{\rho(\lambda,u)d\lambda}{(\lambda_1(u)-\lambda)^2}$. 
Using \eqref{limiting-recursion-relation} , 
we have 
\begin{align}
g^{(4)}(t) &= 6 \int_0^t ds \exp(-2 \int_s^t \varphi ) \int_0^s dr \exp(-\int_r^s \varphi ) h(r) \notag\\
&= 6 \int_0^t ds \exp(- \int_s^t \varphi ) \int_0^s dr \exp(-\int_r^t \varphi ) h(r)  \label{second-line-g-4}
\end{align}
where we have used the Chasles relation in the second line. 
Inverting the order of integration over $s$ and $r$, we get 
\begin{align*}
g^{(4)}(t) = 6 \int_0^t dr \exp(- \int_r^t \varphi ) \int_r^t ds \exp(-\int_s^t \varphi ) h(s)\,. 
\end{align*}
The idea is now to use again the Chasles relation for the second integral to obtain 
\begin{align}\label{g-4-re-appears}
g^{(4)}(t) = 6 \int_0^t dr \exp(- \int_r^t \varphi ) \left(g^{(2)}(t) - \int_0^r ds h(s) \exp(-\int_s^t \varphi )\right) 
\end{align}
which may be rewritten 
\begin{align*}
g^{(4)}(t) = 6 \, g^{(2)}(t)^2 - g^{(4)}(t) 
\end{align*}
where we have noticed from the expression \eqref{second-line-g-4} that the second term of \eqref{g-4-re-appears} is indeed $g^{(4)}(t)$.
The relation $g^{(4)}(t) = 3 g^{(2)}(t)$ follows.

We now have to treat the general case $n\ge 3$ which amounts to prove that 
\begin{align*}
g^{(2n)}(t) = (2n-1) (2n-3) \cdots 3 \, g^{(2)}(t)^n\,. 
\end{align*}
We do the proof recursively on $n$. As a warm up, we first present the proof for $n=3$. 
%The computations are similar as in the $n=2$ case 
%with a few additional tricks, that are already 
%needed to handle for $n=3$.  We only present the proof in the $n=3$ case. 
Using the relation $g^{(4)}(t)= 3 \, g^{(2)}(t)^2$ and \eqref{limiting-recursion-relation}, we have 
\begin{align*}
g^{(6)}(t) &= 45 \int_0^t ds \exp(- 3 \int_s^t \varphi)\,  h(s) \, g^{(2)}(s)^2  \\
&= 45  \int_0^t ds \exp(- 3 \int_s^t \varphi) \int_0^s dr \int_0^s dv \exp(- \int_r^s\varphi )  \exp(- \int_v^s\varphi ) h(r) h(v) \,. 
\end{align*}
Using the Chasles relation as before, we get 
\begin{align}\label{good-form-g-6}
g^{(6)}(t) = 45  \int_0^t ds \exp(-  \int_s^t \varphi) \int_0^s dr \int_0^s dv \exp(- \int_r^t\varphi )  \exp(- \int_v^t\varphi ) h(r) h(v) \,.
\end{align}
Using symmetry properties, we can write 
 \begin{align*}
g^{(6)}(t) = 90  \int_0^t ds \exp(- 3 \int_s^t \varphi) \int_0^s dr \int_0^r dv \exp(- \int_r^s\varphi )  \exp(- \int_v^s\varphi ) h(r) h(v) \,.
\end{align*}
Now we change the order of integration to obtain 
 \begin{align*}
g^{(6)}(t) = 90  \int_0^t dr h(r) \exp(-  \int_r^t \varphi)  \int_0^r dv h(v)  \exp(- \int_v^t\varphi )  \int_r^t ds h(s) \exp(- \int_s^t\varphi )  \,.
\end{align*}
The idea is again to rewrite the last integral thanks to the Chasles relation 
\begin{align}\label{almost-g-6}
g^{(6)}(t) = 90  \int_0^t dr h(r) \exp(-  \int_r^t \varphi)  \int_0^r dv h(v)  \exp(- \int_v^t\varphi ) \left(g^{(2)}(t) - \int_0^r ds h(s) \exp(-\int_s^t\varphi) \right). 
\end{align}
The second integral can be reckoned from the expression \eqref{good-form-g-6} while the first one is easily rewritten in terms of $g^{(4)}(t)$. 
Eq. \eqref{almost-g-6} can thus be rewritten 
\begin{align*}
g^{(6)}(t) = 15 \, g^{(2)}(t) g^{(4)}(t) - 2 g^{(6)}(t)\,,
\end{align*} 
and therefore $g^{(6)}(t) = 15 \, g^{(2)}(t)^3$.

The general case is in fact very similar to the $n=3$ case up to an easy generalization but we write down the proof 
for completeness.  
If $g^{(2n-2)}(t)= (2n-3) (2n-5)\cdots 3 \, g^{(2)}(t)^{n-1}$, then, using \eqref{limiting-recursion-relation},
\begin{align*}
&g^{(2n)}(t) = n (2n-1) (2n-3) \cdots 3  \int_0^t ds\exp\left(- n\int_s^t \int_\R \varphi \right) g^{(2)}(s)^{n-1} \, h(s) \, ds \\
&  = n (2n-1) (2n-3) \cdots 3  \int_0^t ds\exp\left(- n\int_s^t \int_\R \varphi \right)  \int_0^{s} dr_1 h(r_1) \exp(- \int_{r_1}^s \varphi) \\
& \cdots \int_{0}^s dr_{n-1}h(r_{n-1}) \exp(- \int_{r_{n-1}}^s \varphi)\,. 
\end{align*} 
Using the Chasles relation, 
\begin{align}\label{good-form-g-n}
&g^{(2n)}(t) = n (2n-1) (2n-3) \cdots 3  \int_0^t ds\exp\left(- \int_s^t \int_\R \varphi \right)  \int_0^{s} dr_1 h(r_1) \exp(- \int_{r_1}^t \varphi) \notag\\
& \cdots \int_{0}^s dr_{n-1}h(r_{n-1}) \exp(- \int_{r_{n-1}}^t \varphi)\,. 
\end{align}
Using symmetry properties, 
\begin{align*}
&g^{(2n)}(t) =  n \, (n-1) \, (2n-1) (2n-3) \cdots 3  \int_0^t ds \exp\left(- \int_s^t \int_\R \varphi \right)  \int_0^{s} dr_1 h(r_1) \exp(- \int_{r_1}^t \varphi) \\
& \int_0^{r_1} dr_2 h(r_2) \exp(- \int_{r_2}^t \varphi)  \cdots \int_{0}^{r_1} dr_{n-1}h(r_{n-1}) \exp(- \int_{r_{n-1}}^t \varphi)\,. 
\end{align*}
Now we change the order of integration to obtain 
\begin{align*}
&g^{(2n)}(t) = n \, (n-1) \, (2n-1) (2n-3) \cdots 3  \int_0^t 
dr_1 h(r_1)  \exp\left(- \int_{r_1}^t \int_\R \varphi \right)  
 \int_0^{r_1} dr_2 h(r_2) \exp(- \int_{r_2}^t \varphi) \\ & \cdots \int_{0}^{r_1} dr_{n-1}h(r_{n-1}) \exp(- \int_{r_{n-1}}^t \varphi)
 \int_{r_1} ^t ds\,  h(s)\,  \exp(- \int_{s}^t \varphi) \,. 
\end{align*}
The idea is again to rewrite the last integral thanks to the Chasles relation 
\begin{align}\label{almost-g-n}
&g^{(2n)}(t) = n \, (n-1) \, (2n-1) (2n-3) \cdots 3  \int_0^t 
dr_1 h(r_1)  \exp\left(- \int_{r_1}^t \int_\R \varphi \right)  
 \int_0^{r_1} dr_2 h(r_2) \exp(- \int_{r_2}^t \varphi)\notag \\ & \cdots \int_{0}^{r_1} dr_{n-1}h(r_{n-1}) \exp(- \int_{r_{n-1}}^t \varphi)
 \left(g^{(2)}(t) - \int_{0} ^{r_1} ds \, h(s) \, \exp(- \int_s^t \varphi) \right)\,. 
\end{align}
Comparing this expression with \eqref{good-form-g-n} for $n$ and $n-1$, we can rewrite this later equation as 
\begin{align*}
g^{(2n)}(t) = n (2n-1) g^{(2n-2)}(t) g^{(2)}(t) - (n-1) g^{(2n)}(t)\,,
\end{align*}
from which we easily derive the relation $g^{(2n)}(t)= (2n-1) g^{(2)}(t) g^{(2n-2)}(t)$ at level $n$.

\end{document}